\documentclass[leqno,11pt]{amsart}
\pdfoutput=1
\usepackage{fourier}
\usepackage[margin=1in]{geometry}
\usepackage{amssymb}
\usepackage{amsfonts}
\usepackage{amsmath}
\usepackage{fancyhdr}
\usepackage{subfig}
\usepackage{graphicx}
\usepackage{caption}
\usepackage{latexsym}
\usepackage{url}
\usepackage{hyperref}
\usepackage{mathrsfs}
\usepackage{color}
\usepackage{enumitem}
\usepackage{amsthm}
\usepackage{indentfirst}

\theoremstyle{plain}
\newtheorem{theorem}{Theorem}[section]
\newtheorem{lemma}[theorem]{Lemma}

\theoremstyle{definition}

\newtheorem{remark}[theorem]{Remark}
\newtheorem{example}{Example}[section]
\numberwithin{equation}{section}
\allowdisplaybreaks

\newcommand{\rd}{\mathrm{d}}

\newcommand{\eps}{\varepsilon}
\newcommand{\veps}{\varepsilon}

\DeclareMathOperator{\tr}{tr}
\numberwithin{equation}{section}

\begin{document}

\title{A Diabatic Surface Hopping Algorithm based on Time Dependent Perturbation Theory and Semiclassical Analysis} 

\author{Di Fang} \address{ Department of Mathematics, University of
  Wisconsin-Madison, Madison, WI 53706, USA} \email{di@math.wisc.edu}

\author{Jianfeng Lu} \address{Departments of Mathematics, Physics, and
  Chemistry, Duke University, Durham, NC 27708, USA}
\email{jianfeng@math.duke.edu} 

\date{\today}

\thanks{This work is supported in part by NSF grants RNMS-1107444
  (KI-Net), DMS-1454939, DMS-1522184, DMS-1107291 (KI-Net), and the Office of the Vice Chancellor for Research and Graduate Education at the University of Wisconsin-Madison with funding from the Wisconsin Alumni Research Foundation. We thank Shi Jin for helpful discussions. }

\begin{abstract}
  Surface hopping algorithms are popular tools to study dynamics of
  the quantum-classical mixed systems.  In this paper, we propose a
  surface hopping algorithm in diabatic representations, based on time dependent
  perturbation theory and semiclassical analysis. The algorithm can be
  viewed as a Monte Carlo sampling algorithm on the semiclassical path
  space for piecewise deterministic path with stochastic jumps between
  the energy surfaces. The algorithm is validated numerically and it
  shows good performance in both weak coupling and avoided crossing
  regimes.
\end{abstract}

\maketitle

\section{Introduction}

The quantum-classical dynamics is known as an important tool to
investigate non-adiabatic systems. Among them, surface hopping
algorithms, in particular fewest switches surface hopping (FSSH) first
proposed by John Tully \cite%
{tully1971, tully1990, tully1998mixed}, have been applied to many
different practical situations. While its success and improvements
under the adiabatic representation have been well documented in both
chemical and mathematical literatures, its performance in the diabatic
picture seems to be less satisfactory, in which case a full theory or
partial explanation remains as an interesting open problem. In fact,
it was once claimed by Tully himself \cite{tullybook} that the surface
hopping algorithm ``does not appear justified with a diabatic
basis''. Many efforts have been devoted to the improvements of surface
hopping algorithm under the diabatic representation (see, e.g.,
\cite{oleg2, oleg1, plasser2012surface, muller1997surface,
  fabiano2008approximate}). The focus of this paper is to propose a
diabatic surface hopping method using asymptotic derivation from the
Schr\"odinger equation in a diabatic basis. 

A nature question to ask is why to consider the diabatic
representation.  There are various reasons with the three main
motivations listed as follows.  First, as mentioned in
\cite{romero1999towards, dia1}, a diabatic representation is quite
attractive or convenient in terms of simulations, because the diabatic
energies and couplings are smooth functions of the nuclear coordinates
involving multiplicative instead of derivative operators. Hence, it
would make more sense to consider the diabtic representation near the
transition zone where the adiabatic potential energy surface is
ill-defined or the adiabatic coupling becomes very large. Moreover, in
some cases a diabatic basis may arise as a nature choice. For example,
consider the spin transitions, and one could identify the potential
energy surfaces with pure spin character and transitions induced by
the spin-orbit interaction, instead of the mixed spin character in the
adiabatic basis on the contrary. Another example is the photoreaction,
where the molecular system is prepared in a diabatic electronic state,
instead of an adiabatic one, by the laser field \cite%
{muller1997surface, domcke1997theory, kouppel1984multimode}. Third, as
is well-known there are at least two mainstreams of surface hopping
algorithms -- the fewest-switch type (in the spirits of FSSH) and
Landau-Zener type (based on the famous Landau-Zener asymptotics
\cite{landau1932theorie, zener1932non}). Fruitful applications as
there are of both types, a clear connection between the two remains
rather mysterious. Due to the diabatic nature of the original form of
the Landau-Zener formula, a rigorously justified diabatic surface
hopping algorithm may bridge the gap for future study, which also
serves as one of our motivations.

In this paper, we propose a new diabatic surface hopping algorithm in
the same spirit of the FSSH algorithm. Our algorithm is derived based
on semiclassical analysis of the Schr\"{o}dinger equation, adapted
from the adiabatic picture in the recent work of one of us
\cite{lu2016frozen, cailu} to the diabatic setting here. This
algorithm not only agrees with a clear intuition but also is justified
by an asymptotic derivation. Inheriting the merits of the
Schr\"{o}dinger equation based derivation, our algorithm is
time-reversible and also the commonly mentioned drawbacks of surface
hopping, such as frustrated hops and decoherence, do not arise as
issues for the algorithm. The dynamics is described by a system of
ODEs with the hoppings represented by a stochastic jump process. In
simulations, we use a direct Monte Carlo method, and the solution of
the Schr\"{o}dinger equation can be simply recovered by a trajectory
average, \textit{aka} a path integral stochastic representation.

The rest of the paper is organized as follows. Starting with an
introduction of the Schr\"{o}dinger equation in the diabatic
representation that is studied here, we describe an intuitive
understanding of the diabatic surface hopping based on time
pertubation theory, and propose a diabatic surface hopping algorithm
that matches the expansion. In Section 3, the algorithm is tested on
various numerical examples, and the sources of errors are
discussed. Finally, in the last two sections, we derive the algorithm
asymptotically and stochastically, where an asymptotic derivation of
the algorithm is given in Section 3, based on the frozen Gaussian
approximation, \textit{aka} the Herman-Kluk propagator, and a
stochastic interpretation is presented in Section 4, where the
solution of the Schr\"{o}dinger equation is represented by a
stochastic path integral.

\section{The Diabatic Surface Hopping Method}

\subsection{Two level Schr\"{o}dinger equation in a diabatic
  representation}

For simplicity of notations, let us focus in this work on a two-level
Schr\"{o}dinger equation in a diabatic representation:
\begin{equation}
i\veps \partial _{t}%
\begin{pmatrix}
u_{0} \\ 
u_{1}%
\end{pmatrix}%
=-\frac{\veps^2}{2}\Delta _{x}%
\begin{pmatrix}
u_{0} \\ 
u_{1}%
\end{pmatrix}%
+%
\begin{pmatrix}
V_{00} & \delta V_{01} \\ 
\delta V_{10} & V_{11}%
\end{pmatrix}%
\begin{pmatrix}
u_{0} \\ 
u_{1}%
\end{pmatrix}%
,  \label{schd}
\end{equation}%
with initial condition 
\begin{equation}
u(0,x)=u_{\text{in}}(x)=%
\begin{pmatrix}
u_{0}(0,x) \\ 
0%
\end{pmatrix}%
,  \label{initial_data}
\end{equation}%
where $\left( t,x\right) \in \mathbb{R}^{+}\times \mathbb{R}^{m}$,
$u\left( t,x\right) =\left[ u_{0}\left( t,x\right), \ u_{1}\left(
    t,x\right) \right]^{\top}$
is the (two component) wave function, and we have assumed without loss
of generality that the initial condition concentrates on the first
energy surface. The Hamiltonian is given by $H=-\frac{%
  \varepsilon ^{2}}{2}\Delta _{x}+V\left( x\right)$, where
\[
V(x) = 
\begin{pmatrix}
  V_{00}(x) & \delta  V_{01}(x) \\
  \delta V_{10}(x) & V_{11}(x)
\end{pmatrix}
\] 
is a matrix potential. $V$ is Hermitian and hence
$V_{01} = V_{10}^{\ast}$ for every $x$.

In this setup, we have two parameters: $\varepsilon$ is the
semiclassical parameter and $\delta$ is a scaling parameter for the
amplitude of the off-diagonal terms in the matrix potential. A small
$\delta$ means weak coupling between the two diabatic surfaces.  For a
given problem in practice, $\veps$ and $\delta$ are small but fixed
parameters. Here for our algorithm, we would consider the regimes that
$\veps$ and $\delta$ are small, but we will not delve much into the
asymptotic limit that they go to $0$.  The asymptotic limits of these
parameters going to zero are studied for instance for the scaling that
$\delta = \sqrt{\veps} \to 0$ \cite{hagedorn1998landau}, where the
Landau-Zener asymptotic limit was studied; see also \cite{Lasser2005,
  lasser2007, KammererLasser}. Another well known case is the scaling
limit of $\delta \to 0$ with fixed $\veps$, which is related to the
celebrated Marcus theory \cite{marcus1956theory}. We will leave the
rigorous study of the asymptotic limit of our algorithm to future
works; while we will test the algorithm on these cases in our
numerical examples.

\subsection{Surface hopping from the view of time dependent perturbation theory}

As usual for diabatic surface hopping algorithms, we view the diagonal
elements of the potential matrix $V(x)$ as two potential energy
surfaces and the off-diagonal elements as coupling terms. One key
observation behind our diabatic surface hopping is a time dependent
perturbation theory expansion of the solution to the Schr\"odinger
equation. To start with, let us rewrite
\eqref{schd} to separate the diagonal and off-diagonal terms
\begin{equation}
i\varepsilon \partial _{t}%
\begin{pmatrix}
u_{0} \\ 
u_{1}%
\end{pmatrix}%
=%
\begin{pmatrix}
-\frac{\varepsilon ^{2}}{2}\Delta _{x}+V_{00} & 0 \\ 
0 & -\frac{\varepsilon ^{2}}{2}\Delta _{x}+V_{11}%
\end{pmatrix}%
\begin{pmatrix}
u_{0} \\ 
u_{1}%
\end{pmatrix}%
+ \delta
\begin{pmatrix}
  0 &  V_{01} \\
  V_{10} & 0%
\end{pmatrix}%
\begin{pmatrix}
u_{0} \\ 
u_{1}%
\end{pmatrix}.%
\end{equation}%
Then the wave function
$u\left( t,x\right) =\left[ u_{0}\left( t,x\right), u_{1}\left(
    t,x\right) \right]^{\top}$
can be given by the Duhamel's principle by viewing the last term above
as a forcing term
\begin{equation}
u\left( t,x\right) =\mathcal{U}\left( t,0\right) u\left( 0,x\right)
- \frac{i\delta}{\veps} \int_{0}^{t}\mathcal{U}\left( t,t_{1}\right) 
\begin{pmatrix}
0 & V_{01} \\ 
V_{10} & 0%
\end{pmatrix}%
u\left( t_{1},x\right) dt_{1},  \label{duh}
\end{equation}%
where 
\begin{equation*}
\mathcal{U}\left( t,s\right) =%
\begin{pmatrix}
\mathcal{U}_{0}\left( t,s\right)  &  \\ 
& \mathcal{U}_{1}\left( t,s\right) 
\end{pmatrix}%
=%
\begin{pmatrix}
e^{-\frac{i}{\varepsilon }\left( t-s\right) H_{0}} &  \\ 
& e^{-\frac{i}{\varepsilon }\left( t-s\right) H_{1}}%
\end{pmatrix}%
\end{equation*}%
is the propagator of the diagonal part of the equation%
\begin{equation*}
i\varepsilon \partial _{t}%
\begin{pmatrix}
u_{0} \\ 
u_{1}%
\end{pmatrix}%
=%
\begin{pmatrix}
-\frac{\varepsilon ^{2}}{2}\Delta _{x}+V_{00} &  \\ 
 & -\frac{\varepsilon ^{2}}{2}\Delta _{x}+V_{11}%
\end{pmatrix}%
\begin{pmatrix}
u_{0} \\ 
u_{1}%
\end{pmatrix}%
=:%
\begin{pmatrix}
H_{0} &  \\ 
& H_{1}%
\end{pmatrix}%
\begin{pmatrix}
u_{0} \\ 
u_{1}%
\end{pmatrix}%
.
\end{equation*}%
The last equality gives the definition of $H_{k}$, the Hamiltonian of
the $k-$th diabatic energy surface.  To simplify the notations, we
will use below the shorthand for the off-diagonal matrix potential
\begin{equation*}
  M :=
\begin{pmatrix}
0 & V_{01} \\ 
V_{10} & 0%
\end{pmatrix}.%
\end{equation*}%
Substitute \eqref{duh} recursively into the right hand side, we get
\begin{align}
  u\left( t,x\right) & =\mathcal{U}\left( t,0\right) u\left( 0,x\right)
                       - \frac{i\delta}{\veps} \int_{0}^{t}dt_{1}\;\mathcal{U}\left( t,t_{1}\right) M \ \mathcal{U}%
                       \left( t_{1},0\right) u\left( 0,x\right)   \notag   \\
                     & + \biggl(\frac{-i\delta}{\veps}\biggr)^2 \int_{0}^{t}dt_{2}\int_{0}^{t_{2}}dt_{1}\;\mathcal{U}\left(
                       t,t_{2}\right) M \ \mathcal{U}\left( t_{2},t_{1}\right) M\ \mathcal{U}%
                       \left( t_{1},0\right) u\left( 0,x\right)  \label{series} \\
                     & + \biggl(\frac{-i\delta}{\veps}\biggr)^3 \int_{0}^{t}dt_{3}\int_{0}^{t_{3}}dt_{2}\int_{0}^{t_{2}}dt_{1}\;\mathcal{U%
                       }\left( t,t_{3}\right) M\ \mathcal{U}\left( t_{3},t_{2}\right) M\ 
                       \mathcal{U}\left( t_{2},t_{1}\right) M\ \mathcal{U}\left( t_{1},0\right)
                       u\left( 0,x\right)   \notag \\
                     & +\cdots .  \notag
\end{align}%
Recall that we have assumed that the initial datum \eqref{initial_data}
is non-zero only in the first component, the diagonal propagator
$\mathcal{U}$ will keep the non-zero structure, while the action of
$M$ will ``flip'' the non-zero entries in the vector. Following this,
a simple calculation shows that all terms on the right hand side of
\eqref{series} containing even number (including $0$) of $M$ terms would
contribute at time $t$ to the first entry of $u(t, x)$, while the rest
contributes to the second entry, \textit{i.e.},
\begin{align}
u_{0}\left( t,x\right) & =\mathcal{U}_{0}\left( t,0\right) u_{0}\left(
0,x\right)   \notag \\
& + \biggl( \frac{-i\delta}{\veps}\biggr)^2 \int_{0}^{t}dt_{2}\int_{0}^{t_{2}}dt_{1}\mathcal{U}_{0}\left(
t,t_{2}\right)  V_{01} \mathcal{U}_{1}\left(
t_{2},t_{1}\right) V_{10} \mathcal{U}_{0}\left(
t_{1},0\right) u_{0}\left( 0,x\right)   \label{series0} \\
&
+ \biggl(\frac{-i\delta}{\veps}\biggr)^4 \int_{0}^{t}dt_{4}\int_{0}^{t_{4}}dt_{3}\int_{0}^{t_{3}}dt_{2}%
\int_{0}^{t_{2}}dt_{1}\;\mathcal{U}_{0}\left( t,t_{4}\right) 
V_{01} \mathcal{U}_{1}\left( t_{4},t_{3}\right)
V_{10} \mathcal{U}_{0}\left( t_{3},t_{2}\right)
V_{01} \times   \notag \\
& \hspace{14em} \times \mathcal{U}_{1}\left( t_{2},t_{1}\right) V_{10}
\mathcal{U}_{0}\left( t_{1},0\right) u_{0}\left( 0,x\right) \notag \\
& +\cdots ,  \notag\\
\intertext{and,}
u_{1}\left( t,x\right)  &= \biggl(\frac{-i\delta}{\veps}\biggr) \int_{0}^{t}dt_{1}\mathcal{U}_{1}\left(
t,t_{1}\right) V_{10} \mathcal{U}_{0}\left( t_{1},0\right)
u_{0}\left( 0,x\right)   \notag \\
&+ \biggl(\frac{-i\delta}{\veps}\biggr)^3 \int_{0}^{t}dt_{3}\int_{0}^{t_{3}}dt_{2}\int_{0}^{t_{2}}dt_{1}\mathcal{U%
}_{1}\left( t,t_{3}\right)  V_{10}  \mathcal{U}_{0}\left(
t_{3},t_{2}\right)  V_{01}  \mathcal{U}_{1}\left(
t_{2},t_{1}\right) V_{10}  \times \label{series1} \\
& \hspace{14em} \times \mathcal{U}_{0}\left(
t_{1},0\right) u_{0}\left( 0,x\right)   \notag \\
&+\cdots .  \notag
\end{align}%
Note that the sequence $\left\{ t_{j}\right\} _{j=1}^{n}$ corresponds
to the time that the non-zero entries in the wave function are
swapped; in between the consecutive $t_j$'s, the wave function
propagates on one of the diabatic energy surfaces, and thus the
switching can be viewed as a hop from one energy surface to the
other. The series \eqref{series0} and \eqref{series1} thus take into
account the contributions from even number of hops (including no hop)
and odd number of hops, respectively.

To make the above observation more transparent, let us consider the
first several terms of both series.  The first term of \eqref{series0}
$\mathcal{%
  U}_{0}\left( t,0\right) u_{0}\left( 0,x\right) $
represents the contribution of $0$ hops, since it corresponds to wave
packets starting on diabatic surface $0$, following the dynamic on
surface $0$ till the final time $t$. For the first term of
\eqref{series1}, for a fixed $t_1$, the integrand corresponds to wave
packets, initially starting on surface $0$, propagating along the same
potential energy surface from time $0$ to $t_{1}$, "hopping" to
surface $1$ at time $t_{1}$ according to some hopping coefficent
$\frac{i\delta}{\eps}V_{01}$ (to be made precise in the next section),
and then propagating on surface $1$ till the final time
$t$. Similarly, the second term of \eqref{series0} contributes to
hopping twice, with hopping times $t_{1}$ and $t_{2}$, and etc.

Therefore, the total wave function is a summation of the contributions
of all number of hops:
\begin{equation}  \label{u_series}
u\left( t,x\right) =%
\begin{pmatrix}
u_{0}^{\left( 0\right) }\left( t,x\right) +u_{0}^{\left( 2\right) }\left(
t,x\right) +u_{0}^{\left( 4\right) }\left( t,x\right) +\cdots  \\ 
u_{0}^{\left( 1\right) }\left( t,x\right) +u_{0}^{\left( 3\right) }\left(
t,x\right) +u_{0}^{\left( 5\right) }\left( t,x\right) +\cdots 
\end{pmatrix}%
,
\end{equation}%
where $u_{0}^{\left( n\right) }\left( t,x\right) $ denotes the
contribution initiated at surface $0$, with $n$ hoppings by time $t$,
which contains an integral with respect to all possible ordered
hopping times $T_{n:1}:=\left( t_{n},\cdots ,t_{1}\right) $. It is a
good place to point out that this leads to a high dimensional integral
to calculate $u_{0}^{\left( n\right) }$ and thus $u$, not mention the
high dimensional phase space integral already involved in the
semiclassical approximation. To deal with such high dimensional
integrals numerically, it is natural to design a Monte Carlo scheme,
combined with semiclassical approximation for the integrand. This is
exactly the idea behind the surface hopping algorithm in this work.

The series expansion \eqref{series0} converges absolutely for any
fixed $\veps$ and $\delta$ provided $V_{01} \in L_x^{\infty}(\mathbb{R}^m)$, since we only integrate over ordered time
sequence and thus a combinatorial factor in the denominator as the
number of hops (off-diagonal $M$) increases. To be specific, one could easily do a loose estimate on each term in this series: 
\begin{align}
\left\Vert u_{0}^{(n)}\right\Vert_{L_x^2(\mathbb{R}^m)} \leq &
\biggl(\frac{\delta}{\veps}\biggr)^n \int_{0}^{t}dt_{n}\int_{0}^{t_{n}}dt_{n-1}\cdots \int_{0}^{t_{2}}dt_{1}
 \left\Vert V_{01} \right\Vert_{L_x^{\infty}(\mathbb{R}^m)} ^n  
 \left\Vert u_0 \left(0,x \right) \right\Vert_{L_x^2(\mathbb{R}^m)} \label{abs_conv} \\
 =& \biggl(\frac{\delta}{\veps}\biggr)^n  \frac{t^n}{n!} \left\Vert V_{01} \right\Vert_{L_x^{\infty}(\mathbb{R}^m)} ^n  \left\Vert u_0 \left(0,x \right) \right\Vert_{L_x^2(\mathbb{R}^m)} , \notag
\end{align}
which implies the absolute convergence of \eqref{u_series} for any fixed $\veps$ and $\delta$, due to the dominant convergence theorem. 
On the other hand, for larger $\delta$ or smaller $\veps$, we need to
sum over more terms in the series expansion to make sure that the
remainder is below certain threshold. While we do not sum over the
series deterministically term by term in our algorithm, this is still
a good indication that the numerical algorithm would face challenge
when $\delta / \veps$ is large, for example, if we consider fix
$\delta$ and send $\veps$ to $0$. More precisely, in our algorithm, a
large $\delta / \veps$ would lead to a large variance of the estimator
and thus requires increasing number of surface hopping trajectories.
In fact, in the literature, it is well known that the diabatic surface
hopping encounters difficulties when the off-diagonal term is large
compared with the semiclassical parameter, see e.g.,
\cite{neria1993semiclassical}. Thus this observation is perhaps not
limited to our own version of the surface hopping algorithm, but can
be viewed as a generic quantitative difficulty index for the diabatic
surface hopping.

\subsection{A diabatic surface hopping algorithm based on frozen Gaussian approximation}

\label{sh}

The diabatic surface hopping algorithm we propose is based on a
stochastic summation strategy of the series expansion in the last
subsection. The asymptotic justification of the algorithm is
postponed to Section~\ref{asym}.  An overview of the algorithm goes as
follows: We initiate an ensemble of trajectories according to the
initial datum $u\left( 0,x\right) $; for each trajectory, it evolves on
one of the diabatic energy surfaces, and stochastically hops to the
other diabatic surface from time to time; finally, when all the
trajectories are evolved to the final desired time, the wave function
(or more practically the desired observable) is estimated by taking
into account the information from all the trajectories. We will
describe below the two key components of the algorithm:
\begin{enumerate}
\item stochastic dynamics of the trajectories; and
\item estimation of the wave function at the final simulation time. 
\end{enumerate}
The initial sampling of the trajectories is analogous to the adiabatic
case, so we refer the readers to \cite{lu2016improved, lu2016frozen}
for details, which we omit here.

Each trajectory evolves in the extended phase space and is thus
characterized by position $Q_t$, momentum $P_t$, and surface index
$l_t$ to keep track of the current potential energy surface. The
equation of motion of $(Q, P)$ is given by
\begin{align} 
  \frac{\rd}{\rd t}Q & =P,  \label{ode1} \\
  \frac{\rd}{\rd t}P & =-\nabla V_{l_{t}l_{t}}(Q) \notag,
\end{align}
where the force depends on the current surface $l_t$. The surface
index $l_{t}$ follows a Markov jump process on the state space
$\left\{0,1\right\}$, where infinitesimal transition probability is given
by%
\begin{equation}
P\left( l_{t+\delta t}=k \ \lvert \ l_{t}=j,Q_{t}=q\right) =\delta _{jk}+\lambda
_{jk}\left( q\right) \delta t+o\left( \delta t\right) ,
\end{equation}%
with transition rate
\begin{equation}
\lambda \left( q\right) =%
\begin{pmatrix}
\lambda_{00}\left( q\right)  & \lambda _{01}\left( q\right)  \\ 
\lambda_{10}\left( q\right)  & \lambda _{11}\left( q\right) 
\end{pmatrix}%
=%
\begin{pmatrix}
  -\frac{\delta}{\eps} \left\rvert V_{01}\left( q\right) \right\rvert & \frac{\delta}{\eps} \left\rvert
    V_{01}\left(
      q\right) \right\rvert  \\
  \frac{\delta}{\eps} \left\rvert V_{10}\left( q\right) \right\rvert & -\frac{\delta}{\eps} \left\rvert
    V_{10}\left( q\right) \right\rvert
\end{pmatrix}%
.  \label{gen}
\end{equation}%
Here the modulus of $V_{01}$ is taken since we interpret that as a
non-negative rate. Therefore, $(Q_{t},P_{t},l_{t})$ can be viewed as a
Markov switching process whose probability distribution
$F_{t}\left(q,p,l\right) $ follows the following Kolmogorov forward
equation
\begin{equation}
\partial _{t}F_{t}\left(q,p,l\right) + p \cdot \nabla _{q}F_{t}\left(q,p,l\right) -\nabla _{p}V_{ll}(q) \cdot \nabla _{p}F_{t}\left(q,p,l\right)=\sum_{j=0}^{l}\lambda _{jl}\left( q\right) F_{t}\left(
q,p,j\right),
\end{equation}%
where the last two terms of the left hand side is induced by the
Hamiltonian flow of $(Q_{t},P_{t})$, and the right hand side stands
for the contribution of the Markov jump process. Note that $l_t$ is a
piecewise constant and almost surely each trajectory has finite number
of jumps. Consider a realization of a trajectory with $n$ jumps, we
will denote the discontinuity set of $l_{t}$ as
$\left\{ t_{j}\right\} _{j=1}^{n}$, which is the ordered set of time
of jumps of the trajectory.

Along each trajectory, we also solve for the classical action $S_t$
and weighting factor $A_t$ given by
\begin{align}  \label{ode2}
  \frac{\rd}{\rd t}S & =\frac{1}{2}\lvert P\rvert ^{2}-V_{l_{t}l_{t}}\left( Q\right) ,   \\
  \frac{\rd}{\rd t}A & =\frac{1}{2}A\, \tr\left(
                  Z^{-1} \left( \partial _{z}P-i\partial _{z}Q\nabla _{Q}^{2}V_{l_{t}l_{t}}\left(
                  Q\right) \right) \right) ,  \notag
\end{align}
where we have used the short hand
 \begin{equation*} 
   \partial _{z}:=\partial _{q}-i\partial _{p}, \quad \text{and} \quad
   Z:=\partial _{z} ( Q + i P).
\end{equation*}
For the trajectory with initial motion
$\left( Q_0, P_0 \right)= \left( q_0, p_0 \right)$, the initial action
 and  amplitude are assigned as
\begin{equation}  \label{initial_A}
S_0 = 0 \qquad \text{and} \qquad 
A_{0} 
=2^{m/2}\int_{\mathbb{R} ^{m}}u_{0}\left(0,
y\right) e^{\frac{i}{\varepsilon }\left( -p_0\left( y-q_0\right) +\frac{i}{2}%
\lvert y-q_0\rvert ^{2}\right) }dy, 
\end{equation}
where the integral for $A_0$ is either analytically given or
numerically estimated.

Finally, the wave function is reconstructed via a trajectory average, that is, an expectation over all of the trajectories,
\begin{equation}  \label{avg1}
\tilde{u}\left( t,x\right) =C_{\mathcal{N}} \mathbb{E}\left[ 
 \left(-i  \right)^n 
\left( \prod_{j=1}^{n}\frac{V_{l_{t_j}l_{t_{j-1}}} \left(Q_{t_j} \right) }{%
\left\lvert V_{l_{t_j}l_{t_{j-1}}} \left(Q_{t_j} \right) \right\rvert }%
\right) \frac{A_{t}}{\lvert
A_{0} \rvert } \exp \left( \frac{i}{%
\varepsilon }\Theta _{t} \right) 
\exp \left( \int_{0}^{t} \frac{\delta}{\eps} \lvert
V_{01} \left(Q_s \right) \rvert ds\right) 
\begin{pmatrix}
\mathbb{I}_{n\ \text{is even}} \\ 
\mathbb{I}_{n\ \text{is odd}}%
\end{pmatrix}
\right] ,
\end{equation}
where the normalization constant reads 
\begin{equation}  \label{cn}
C_{\mathcal{N}}=\frac{1}{\left( 2\pi \varepsilon \right) ^{3m/2}}\int_{%
\mathbb{R}^{2m}}\lvert A_{0}\left(q_{0},p_{0}\right)
\rvert dq_{0}dp_{0},
\end{equation}
and we have used $\mathbb{I}$ to denote the indicator function and
\begin{equation}\label{theta} 
  \Theta_{t} = S_{t} +P_{t} \left(
    x-Q_{t} \right) +\frac{i}{2}\lvert x-Q_{t} \rvert ^{2}.
\end{equation}
While the expression \eqref{avg1} appears complicated, all the terms
involved in the expectation are obtained along the trajectories. In
the algorithm, the expectation over path ensemble will be estimated by
Monte Carlo algorithm sampling the path space, as will be explained in
more details below.

\section{Numerical Implementation and Examples} \label{numerical ex}

\subsection{Numerical implementation}

Numerical description of our scheme is as follows:
Given initial datum $u_{\text{in}}(x)=\left[ u_{0}\left( 0,x\right) \ 0 \right] ^{T}$,
one first calculates the initial amplitude according to \eqref{initial_A} by numerical quadratures (or even analytically, for special initial data). One samples the initial position and momentum according to $\lvert
A_{0}\left( q_{0},p_{0}\right) \rvert $, viewed as a joint density
function of $\left( q_{0},p_{0}\right)$. Each pair of $\left( q_{0},p_{0}\right)$ corresponds to a
trajectory, with the initial position, momentum and action given as 
\begin{equation}
Q_{0} =q_{0},  \ \
P_{0} =p_{0},  \ \
S_{0} =0.  \label{initial_pqs}
\end{equation}
Then one marches the ODE system \eqref{ode1} and \eqref{ode2} by a ODE
integrator (high-order ones are preferred and we choose here the
fourth order Runge-Kutta scheme). For each time step $t=t^{k}$, the
hopping probability of this time step is calculated as $%
\Delta t \delta \lvert V_{01} \left(Q_{t^k} \right) \rvert $,
where $\Delta t$ is the time step size chosen sufficiently small. We
thus draw a random number uniformly distributed on the interval
$[0,1]$; a hop happens if the probability is larger than the drawn
random number, upon which we switch the potential energy surface index
$\ell_t$, with all of the trajectory coefficients $(Q_t,P_t,S_t,A_t)$
being continuous. Moreover, at each of the hopping time $t_{n}$, we
record the hopping coefficients $V_{01}\left(Q_{t_n} \right)$. Repeat
this procedure until the desired final time, the wave function (or
expectation of an observable) is estimated according to the trajectory
average \eqref{avg1}.

\begin{figure}[htbp]
{}
\setcaptionwidth{4.6in}
\begin{center}
\includegraphics[scale=0.6]
{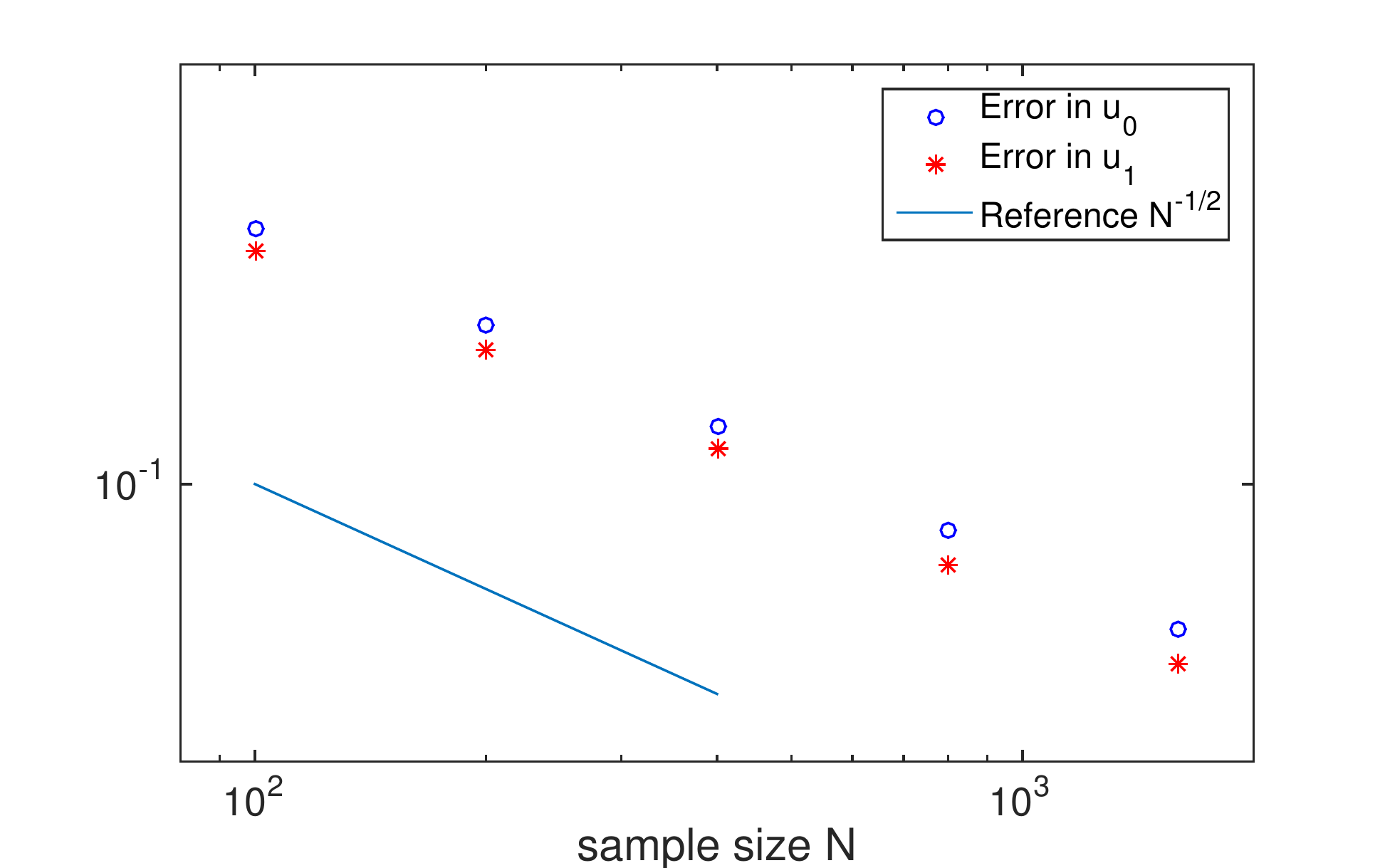}
\end{center}
\caption{Sampling error incurred by the Monte Carlo method. Here
  $\protect%
  \varepsilon$
  is fixed as 0.04, and test on various numbers of trajectories $%
  N=100, 200, 400, 800, 1600$.
  The plot confirms the convergence rate of $%
  N^{-1/2} $
  and that the sampling error dominates the total numerical error.}
\label{conv}
\end{figure}

Before we turn to numerical examples, let us discuss the possible
source of errors of the numerical solution compared to the exact
solution to the Schr\"{o}dinger equation:
\begin{enumerate}[label=(\roman*)]
\item Initial error, coming from the numerical estimation of $A_{0}$,
  if analytic expression is not available. Note that if the phase
  space is high dimensional, numerical quadrature to obtain $A_0$ is
  not feasible, and hence we would need to restrict to specific
  initial data such that analytic expression is available.

\item Asymptotic error, introduced by the semiclassical approximation
  of the evolution on single diabatic surfaces, where we neglect 
  terms with orders higher than $\veps$ in the frozen Gaussian
  approximation. This error is of order
  $\exp{ \left( \frac{\delta T}{\eps} \right) \eps}$, where $T$ is the
  final time (see Section \ref{asymptotic derivation} for further
  discussion). 

\item Numerical integration error of the ODEs. This error can be
  controlled by choosing high order ODE solvers. In our case the
  fourth order Runge-Kutta method seems to provide good numerical accuracy.

\item Sampling error, incurred by the application of the Monte Carlo
  method.  This error decays as $N^{-1/2}$, where $N$ denotes the
  number of trajectories. A numerical validation of the sampling error
  is given in Fig. \ref{conv}. Typically the sampling error dominates
  the total error of the numerical algorithm. 
\end{enumerate}

\subsection{Numerical examples}

In this section, the numerical examples are chosen in the spirit of
the three standard test cases first considered in Tully's paper
\cite{tully1990}. The reference solutions are all computed by a
time-splitting pseudo-spectral scheme \cite{Bao:2002fy} (this is
possible as these test cases are in low spatial dimension), with the
meshing strategy $\Delta t = o(\varepsilon)$,
$\Delta x = o(\varepsilon)$ and refined enough to guarantee that the
reference solution is close to the true solution.

\begin{example}[Simple Crossing] \label{ex_simple} Consider the matrix
  potential %
\begin{equation}  \label{V_simple}
\begin{pmatrix}
\tanh x & \delta \\ 
\delta & -\tanh x%
\end{pmatrix}%
,
\end{equation}%
whose diabatic potential energy surfaces is shown in Fig.~\ref{simple1}. In our test, we choose 
$\varepsilon =\delta =0.04$, the final time $T=1.2$, and the initial datum given by  
\begin{equation*}
u_{0}\left(0, x\right) =e^{-12.5\left( x+1.5\right) ^{2}}e^{\frac{2i}{%
\varepsilon }\left( x+1.5\right) }\text{.}
\end{equation*}%
The number of trajectories $N$ is $5000$. It can be seen in
Fig.~\ref{simple2} that the numerical solutions match nicely with the
reference solutions. This validates that our algorithm can correctly
reconstruct the wave function.

\begin{figure}[htbp]
{}
\begin{center}
\subfloat[Diabatic Potential Energy Surfaces]{\includegraphics[scale=0.4]
{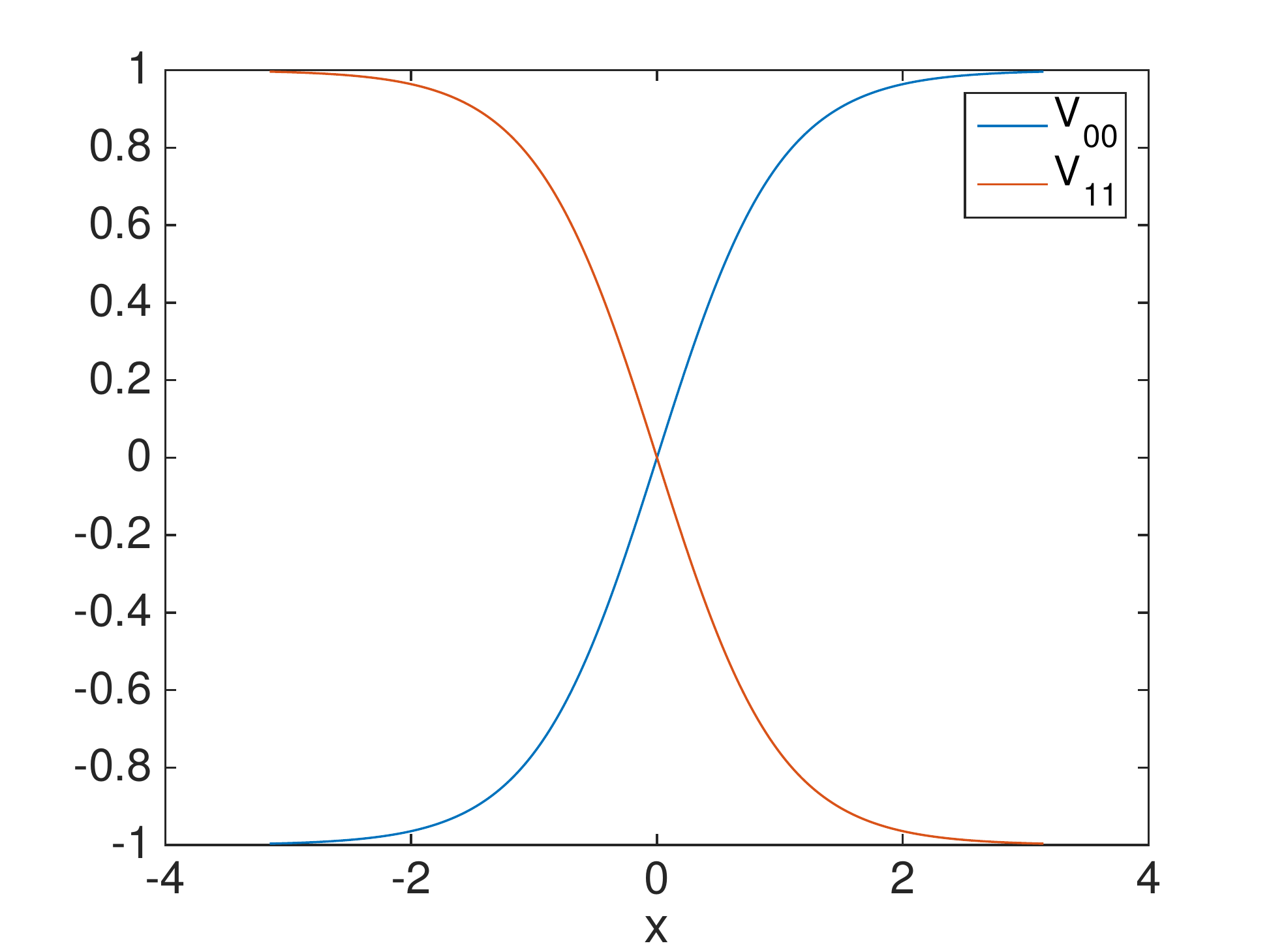} \label{simple1}} 
\subfloat[Numerical solutions]{\includegraphics[scale=0.4]
{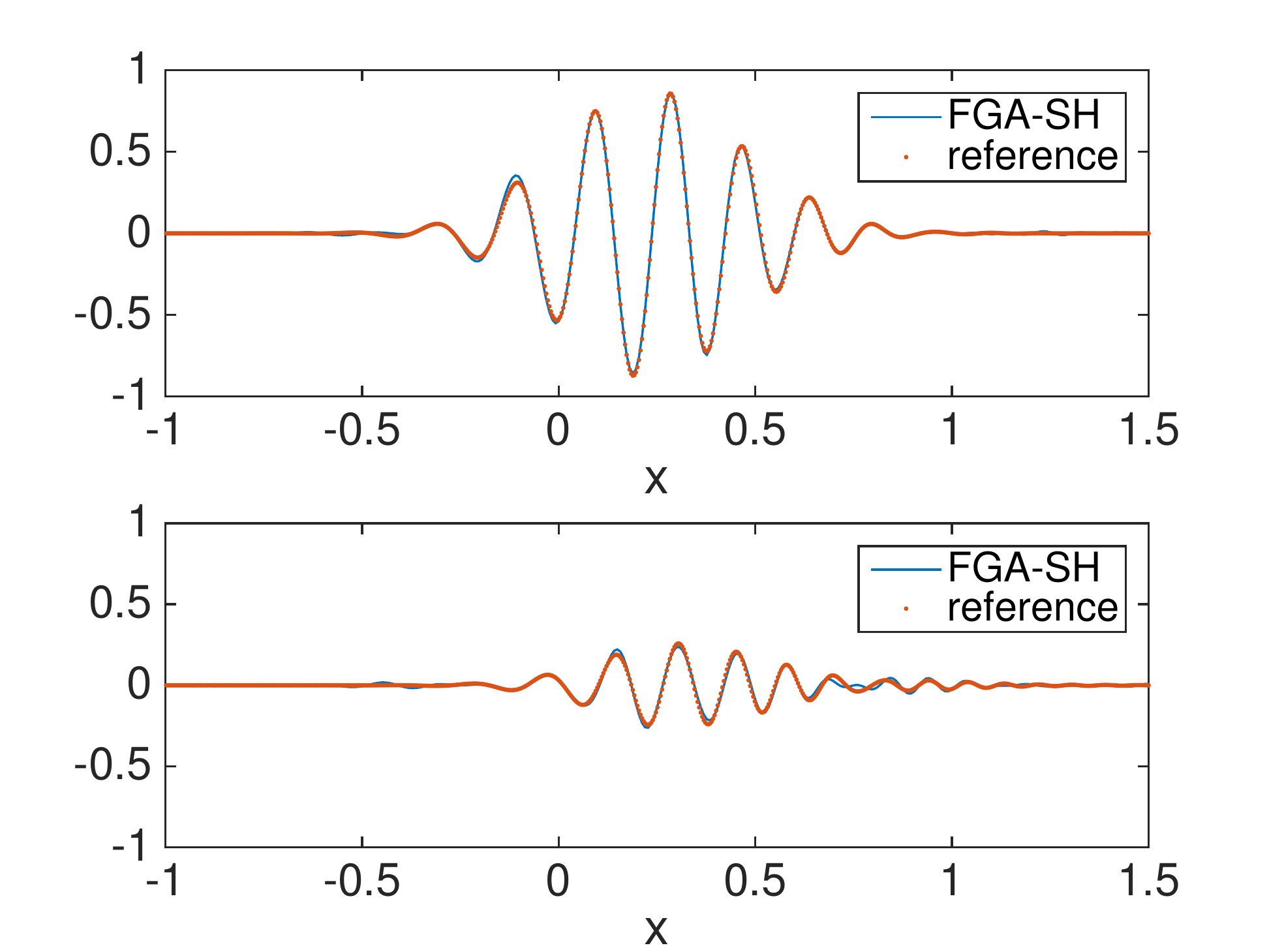} \label{simple2}}
\end{center}
\caption{Simple Crossing. (Left panel): the diabatic potential energy
  surfaces. (Right panel): the numerical solutions (blue solid) of
  $u_{0}$ (top) and $u_{1}$ (bottom) agree well with the reference
  solutions (red dashed). The reference solutions is computed by the
  time-splitting pseudo-spectral scheme with
  $\Delta t = \frac{\protect\varepsilon}{32}$ and
  $\Delta x=\frac{(b-a)\protect\varepsilon}{32}$.}
\label{simple}
\end{figure}
\end{example}

\begin{example}[Dual Crossing]
  In the second example, we consider the case of dual crossing; it is
  a more challenging example as it exhibits two avoided crossing. So
  to get the correct solution, the algorithm has to capture the phase
  of the wave function during the evolution. The matrix potential
  reads%
\begin{equation*}
\begin{pmatrix}
0 & 0.015e^{-0.06x^{2}} \\ 
0.015e^{-0.06x^{2}} & -0.1e^{-0.28x^{2}}+0.05%
\end{pmatrix}%
.
\end{equation*}%
We take in this example $\varepsilon =\frac{1}{\sqrt{2000}}$, 
and initial datum
\begin{equation*}
u_{0}\left(0, x\right) =e^{-\sqrt{500}\left( x+2.5\right) ^{2}}e^{\frac{2i%
}{\varepsilon }\left( x+2.5\right) }.
\end{equation*}%
We run the algorithm until the wave packets pass both crossing areas
around $T=2.2$.  The number of the trajectories is $10000$. For this
example, we observe again that the numerical and reference solutions
are in good agreement.

\begin{figure}[htbp]
\begin{center}
\subfloat[Diabatic Potential Energy Surface]{\includegraphics[scale=0.4]
{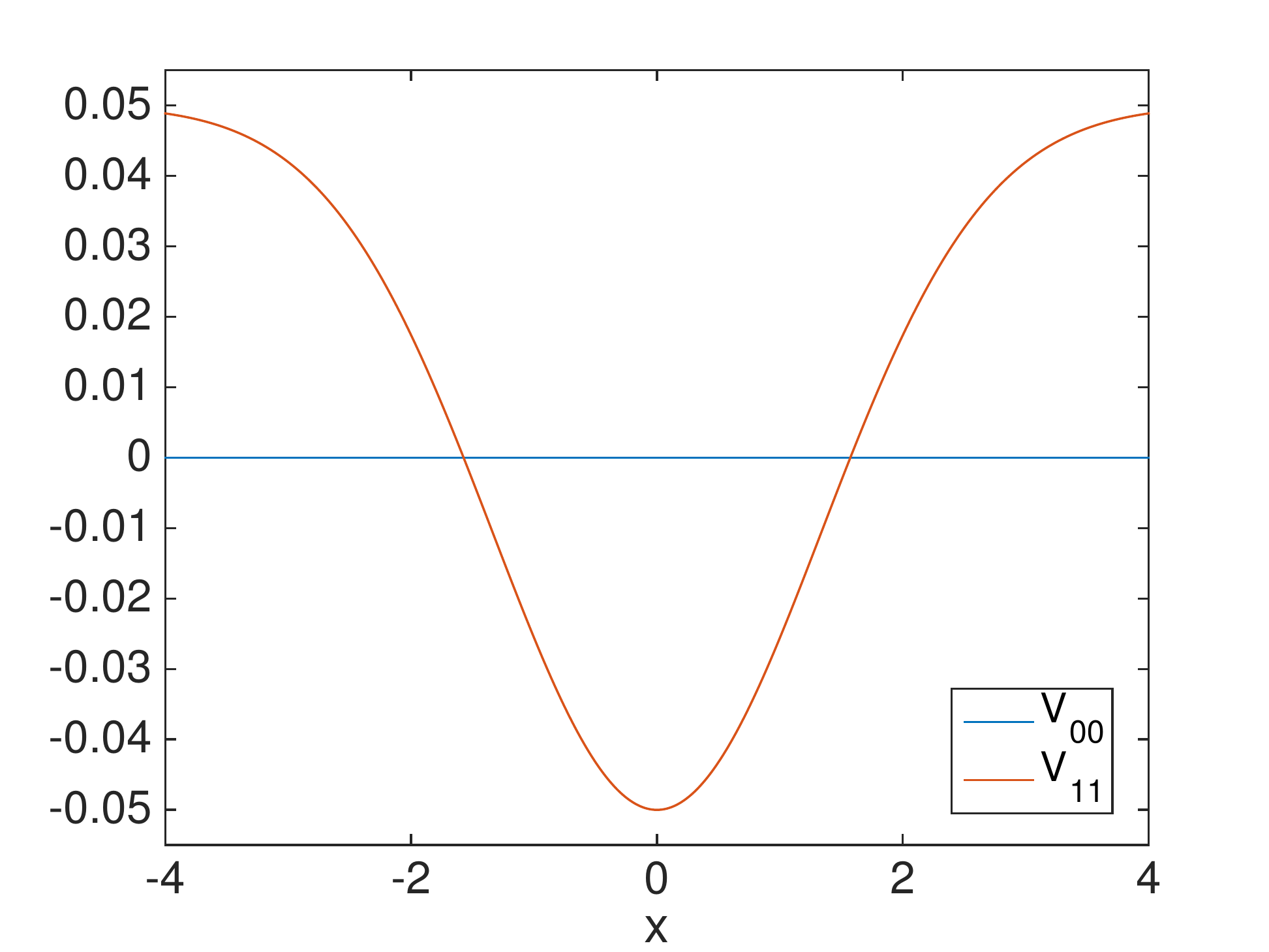}  } 
\subfloat[Numerical solutions]{\includegraphics[scale=0.4]
{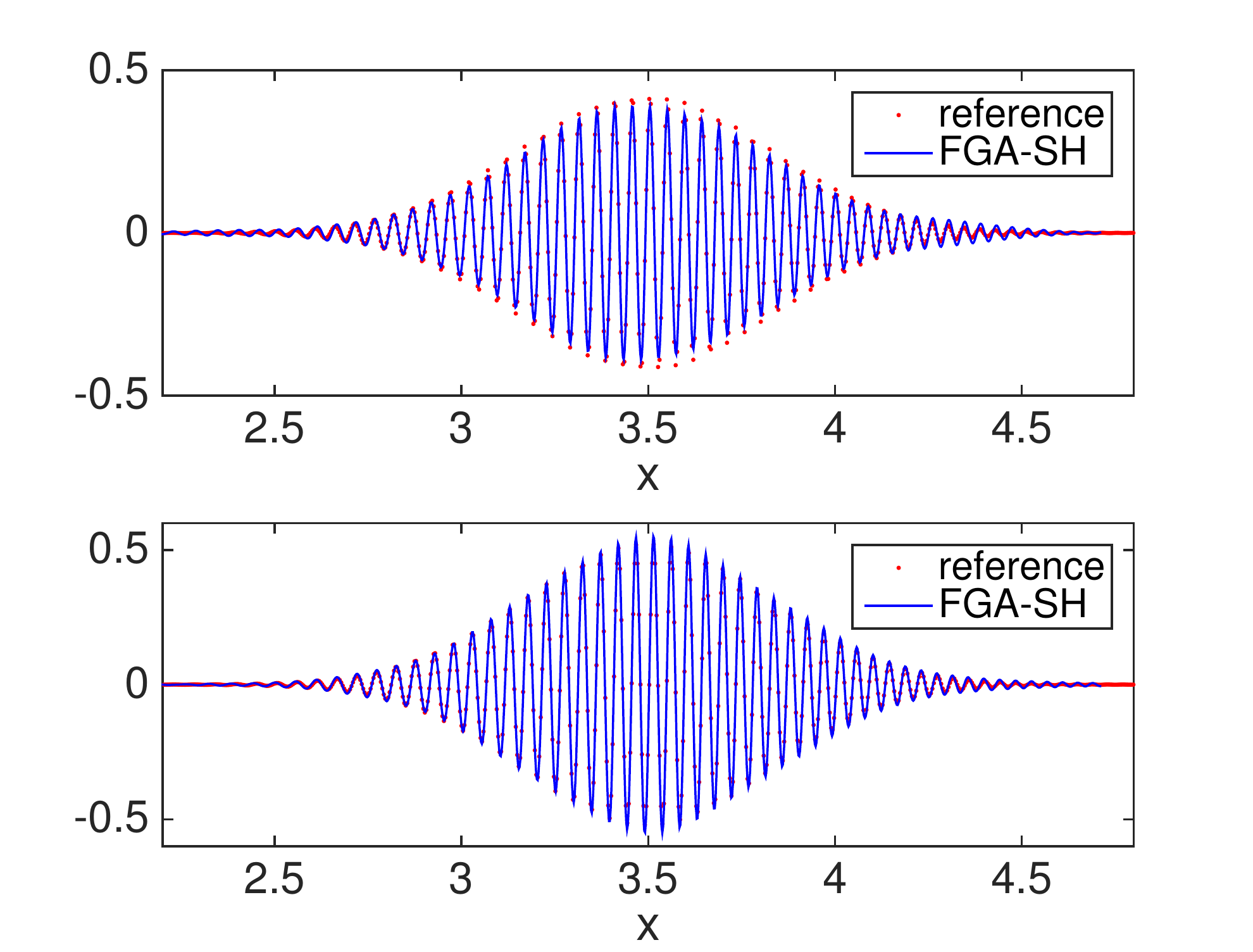}}
\end{center}
\caption{Dual Crossing. (Left panel): the diabatic potential energy
  surfaces. (Right panel): the numerical solutions (blue solid) of
  $u_{0}$ (top) and $u_{1}$ (bottom) agree well with the reference
  solutions (red dashed). The reference solutions is computed by the
  time-splitting pseudo-spectral scheme with
  $\Delta t = \frac{\protect\varepsilon}{32}$ and
  $\Delta x=\frac{(b-a)\protect\varepsilon}{32}$.}
\label{dual}
\end{figure}
\end{example}

\begin{example}[Extended Coupling with Reflection]
  In this example, we consider the matrix potential
\begin{equation*}
\begin{pmatrix}
\arctan \left( 10x\right) +\frac{\pi }{2} & \delta \\ 
\delta & -\arctan \left( 10x\right) -\frac{\pi }{2}%
\end{pmatrix}%
.
\end{equation*}%
This is often considered as a challenging example, since it contains
an extended region of strong nonadiabatic coupling.  An interesting
feature of this example is that the wave packet evolves on the upper
surface would be partially reflected (bounced back), as shown in
Fig. \ref{extended2} (note the two wave packets on the upper diabatic
energy surface). Consider $\varepsilon = \delta =0.04$, the final time
$T=1.4$, and the initial data same as in Example~\ref{ex_simple}. The
number of the trajectories is $30000$. It can be seen in
Fig.~\ref{extended} that our algorithm still keeps great agreement
with the reference solutions. Both the transmitted and reflected wave
packets on the second energy surface are accurately captured.

\begin{figure}[htbp]
{}
\begin{center}
\subfloat[Diabatic Potential Energy Surface]{\includegraphics[scale=0.4]
{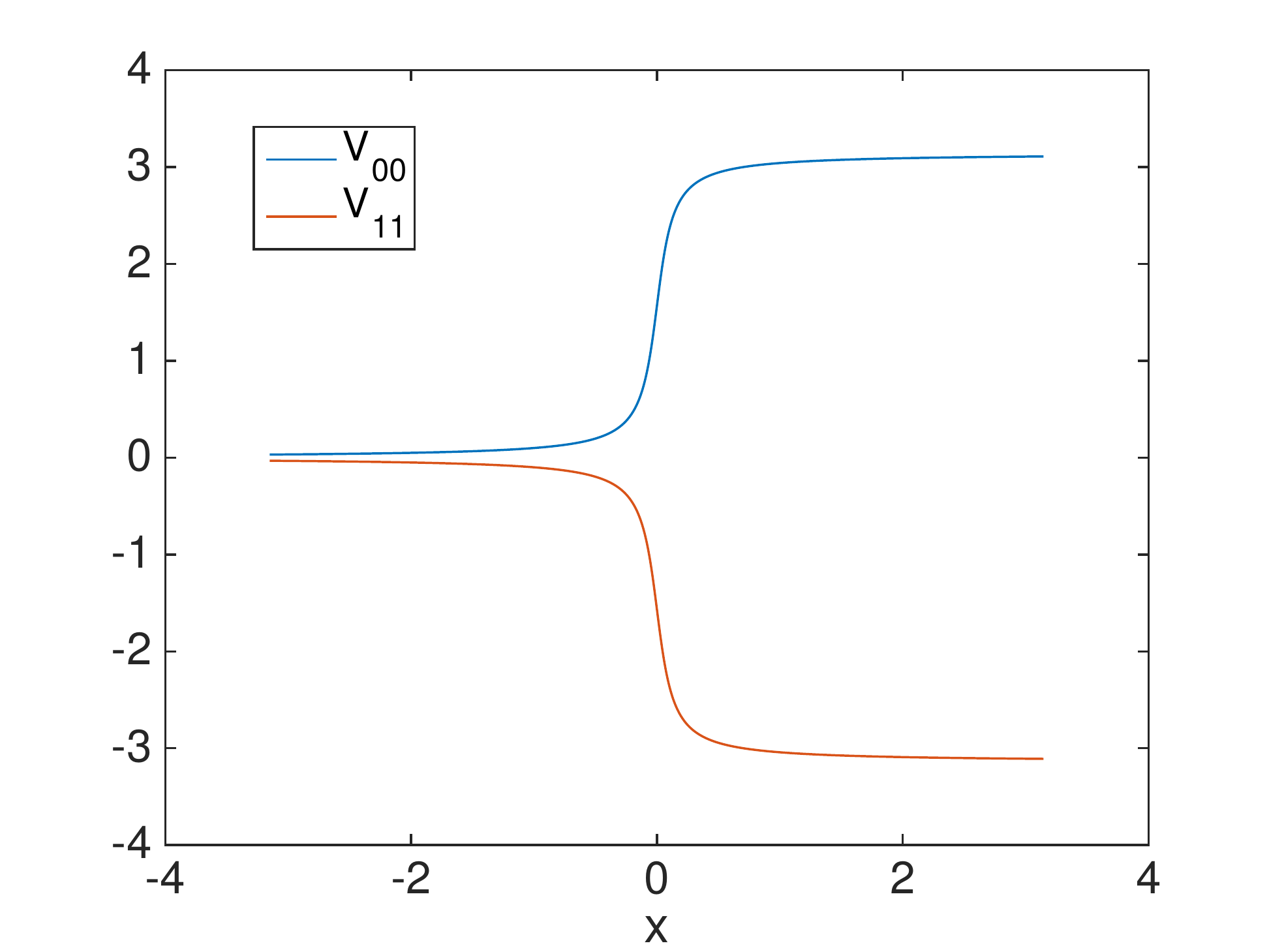}} 
\subfloat[Numerical solutions]{\includegraphics[scale=0.4]
{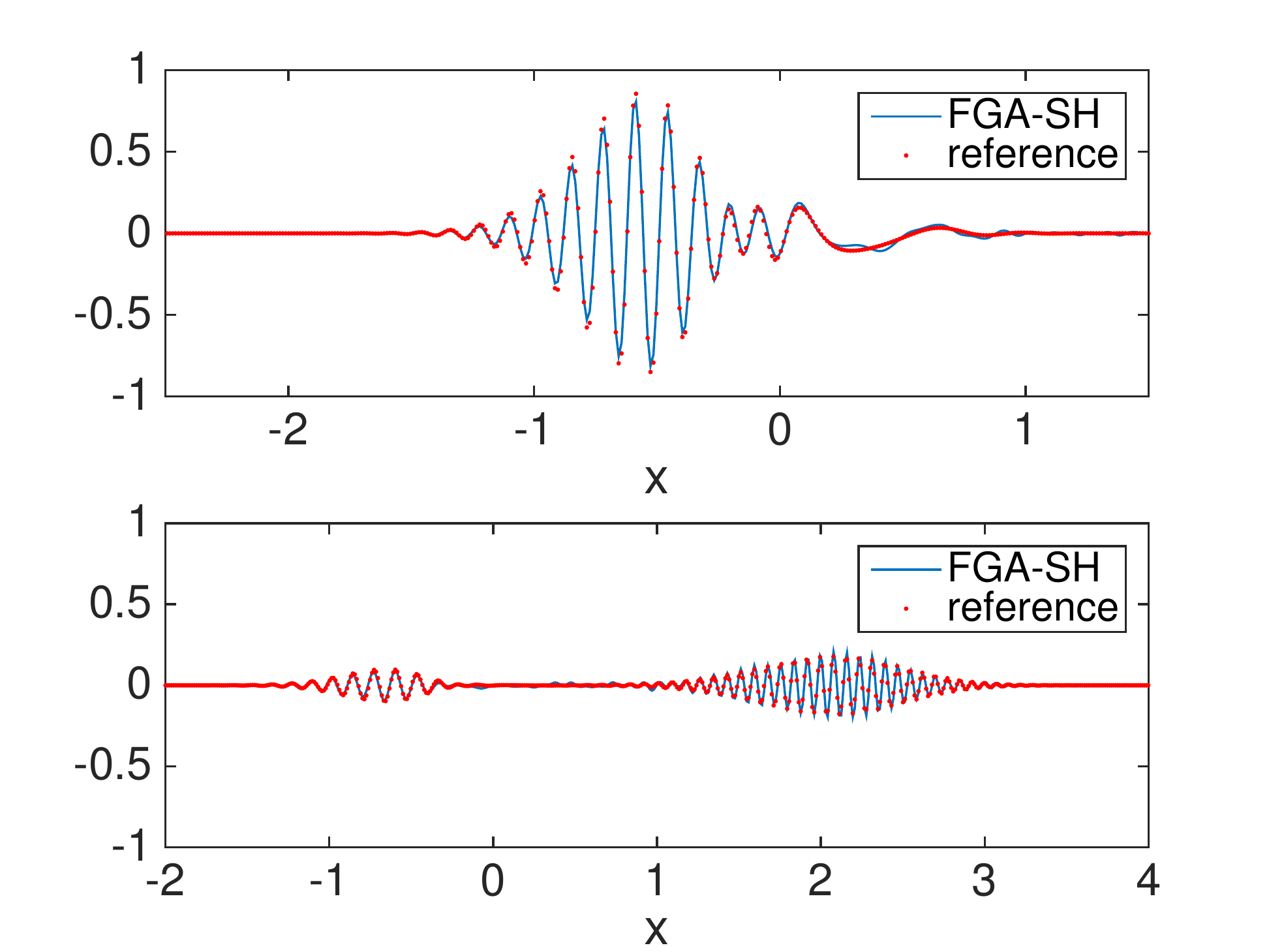} \label{extended2}}
\end{center}
\caption{Extended Coupling with Reflection. (Right panel): the numerical
solutions (blue solid) of
  $u_{0}$ (top) and $u_{1}$ (bottom) show good agreements with the reference solutions
(red dashed). }
\label{extended}
\end{figure}
\end{example}

\begin{example}[Weak Coupling Regime] \label{ex_marcus} In this
  example, we focus on the weak coupling regime when
  $\delta \rightarrow 0$ and $\eps$ fixed. This regime is related to
  the Marcus theory \cite{marcus1956theory} in the theoretical
  chemistry literature, it is particularly interesting since as shown
  in \cite{landry2011communication}, the original Tully's FSSH does
  not capture the right asymptotic limit.  Let us define the
  transition rate as
  \begin{equation}
    R=  \frac{ \int_{\mathbb{R}^{m}} \left \vert u_{1} (T,x) \right \vert ^2 dx}{ \int_{\mathbb{R}^{m}} \left \vert u_{0} (0,x) \right \vert ^2 dx},
  \end{equation}
  where $T$ is some time when the wave packets bypass the transition
  zone. 
%
 The Marcus theory states that the transition rate is
  proportional to the square of the off diagonal term
  $\delta \lvert V_{01} \rvert$. Unfortunately, as is shown in
  \cite{landry2011communication}, the standard FSSH gives the rate
  proportional to $\delta \lvert V_{01} \rvert$ instead of
  $\delta^2 \lvert V_{01} \rvert^2$. Our algorithm is able to capture the
  correct asymptotic limit, as shown in Fig.~\ref{marcus}; the
  numerically predicted transition rate is quadratic in
  $\delta \left\lvert V_{01} \right\rvert$. In this test, we take the matrix
  potential \eqref{V_simple} with various
  $\delta$ 
  and initial datum
  \begin{equation*}
    u_{0}\left(0, x\right) =e^{-12.5\left( x+1\right) ^{2}}e^{\frac{2i}{%
        \varepsilon }\left( x+1\right) }
  \end{equation*}
  for the final time $T=1$.

\begin{figure}[htbp]
{}
\begin{center}
\includegraphics[scale=0.5]
{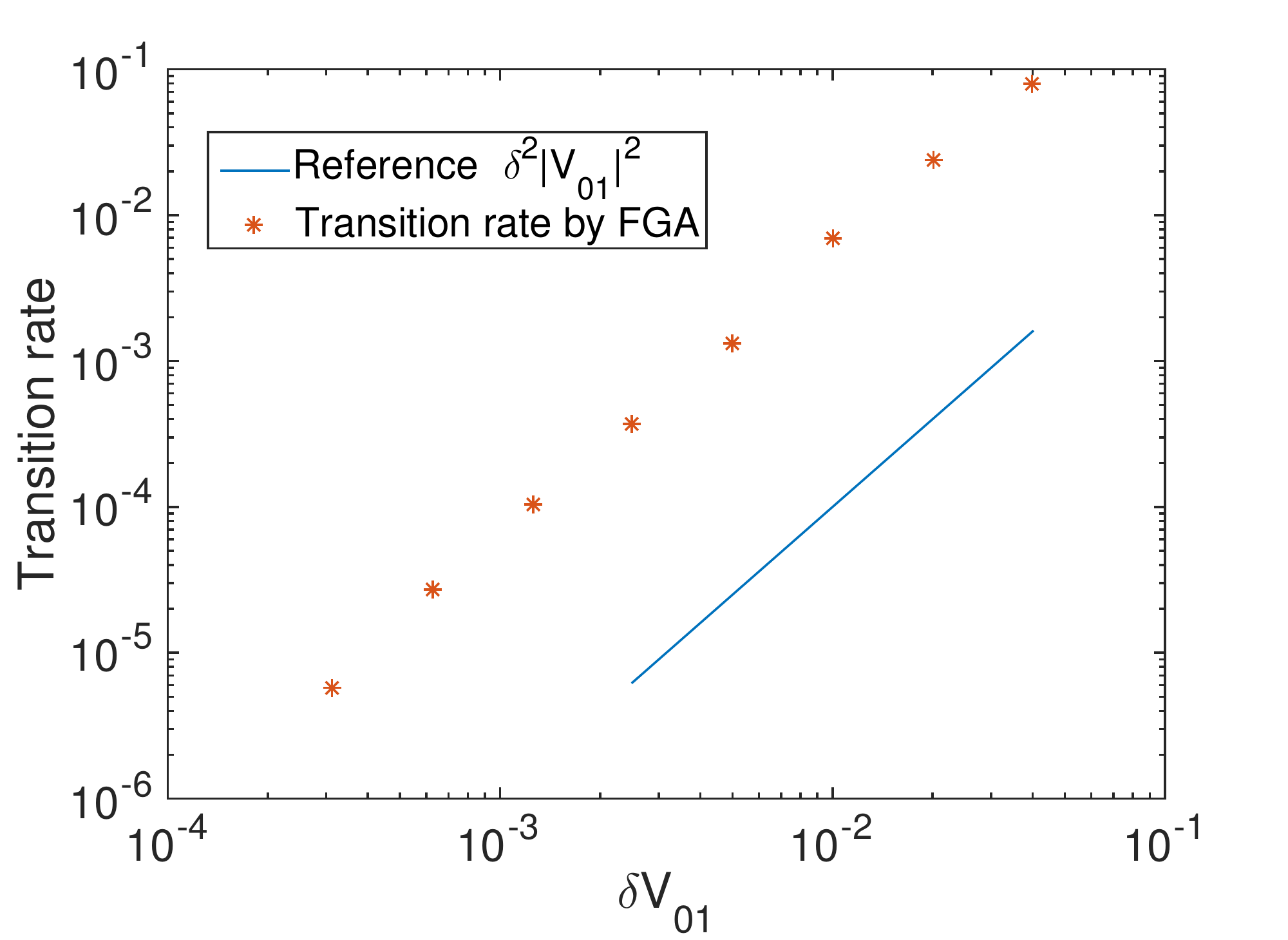}
\end{center}
\caption{Weak Coupling Regime. The transition rate computed by our numerical
  scheme (red~$\ast$) grows proportional to $\delta^2 \lvert V_{01} \rvert^2$,
  which agrees with the Marcus theory. }
 \label{marcus}
\end{figure}

\end{example}

\begin{example}[Landau-Zener Regime] 
  We shall now test the performance of our algorithm in the
  Landau-Zener regime as in \cite{hagedorn1998landau}, i.e., when the
  off-diagonal term of the potential matrix $V$ is
  $O(\eps^{\frac{1}{2}})$. Consider the matrix potential
\begin{equation} 
\begin{pmatrix}
\tanh x & \varepsilon^{\frac{1}{2}} \\ 
\varepsilon^{\frac{1}{2}} & -\tanh x%
\end{pmatrix},%
\end{equation}
with the semiclassical parameter $\eps=0.04$, final time $T=0.5$, and
initial datum given by
\begin{equation*}  \label{initial_p3}
u_{0}\left(0, x\right)=e^{-12.5\left( x+1\right) ^{2}}e^{\frac{3i%
}{\varepsilon }\left( x+1\right) }.
\end{equation*}
We take $20000$ trajectories in the numerical test. As shown in
Fig.~\ref{lz_case}, our scheme still produces good agreement with the
reference solution. On the other hand, compared to
Example~\ref{ex_simple}, where the off-diagonal terms of the matrix
potential has magnitude $O(\eps)$, we do need significantly more
trajectories to achieve similar accuracy. We thus further study in the
next two examples the relation of the number of trajectories versus
the ratio $\delta / \veps$.

\begin{figure}[htbp]
{}
\begin{center}
\includegraphics[scale=0.4]
{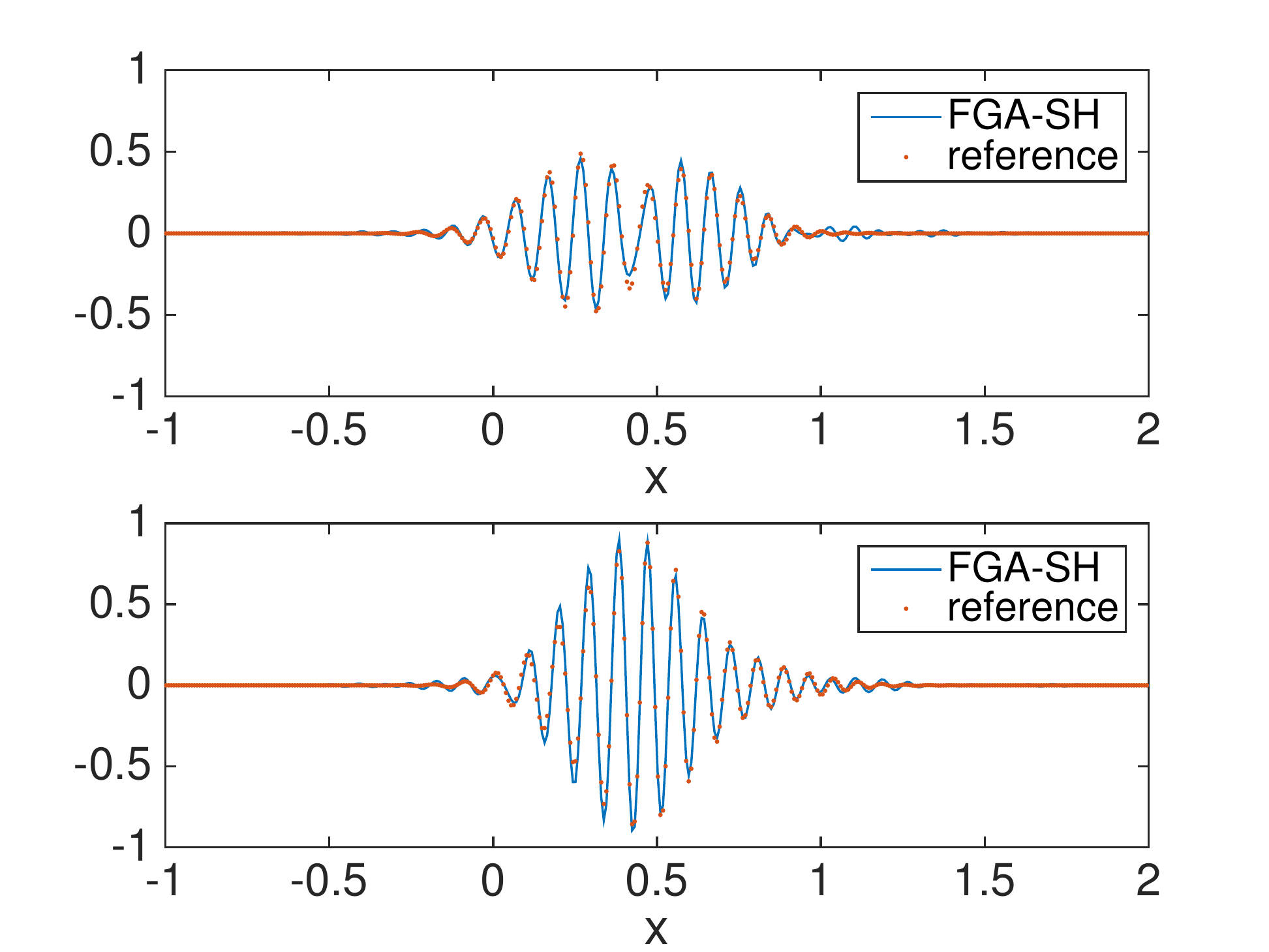} 
\end{center}
\caption{Landau-Zener Regime. The numerical
solutions (blue solid) show good agreements with the reference solutions
(red dashed). }
\label{lz_case}
\end{figure}
\end{example}

\begin{example}[Number of Trajectories versus $\frac{\delta}{\eps}$]
  As have seen in the previous example, more trajectories are needed
  as $\frac{\delta}{\eps}$ becomes larger. In this example, we
  consider fixing $\eps=0.02$ and increasing the values of
  $\delta$. The potential used in this example is \eqref{V_simple}
  (the previous example corresponds to $\delta = \veps^{1/2}$) and
  initial datum
  \begin{equation*}  \label{initial_p3}
u_{0}\left(0, x\right)=e^{-25\left( x+.4\right) ^{2}}e^{\frac{2i%
}{\varepsilon }\left( x+.4\right) }
\end{equation*}
for the time $T=0.5$. We record the number of
  trajectories needed for the $L^2$ error of the wave function to
  achieve certain fixed error threshold. One can see that the number of trajectories grows with respect to $\frac{\delta}{\eps}$ exponentially. As predicted before, it becomes more challenging as $\frac{\delta}{\eps}$ become larger, and one would need more trajectories to get accurate wave function. 
\begin{figure}[htbp]
{}
\begin{center}
\includegraphics[scale=0.5]
{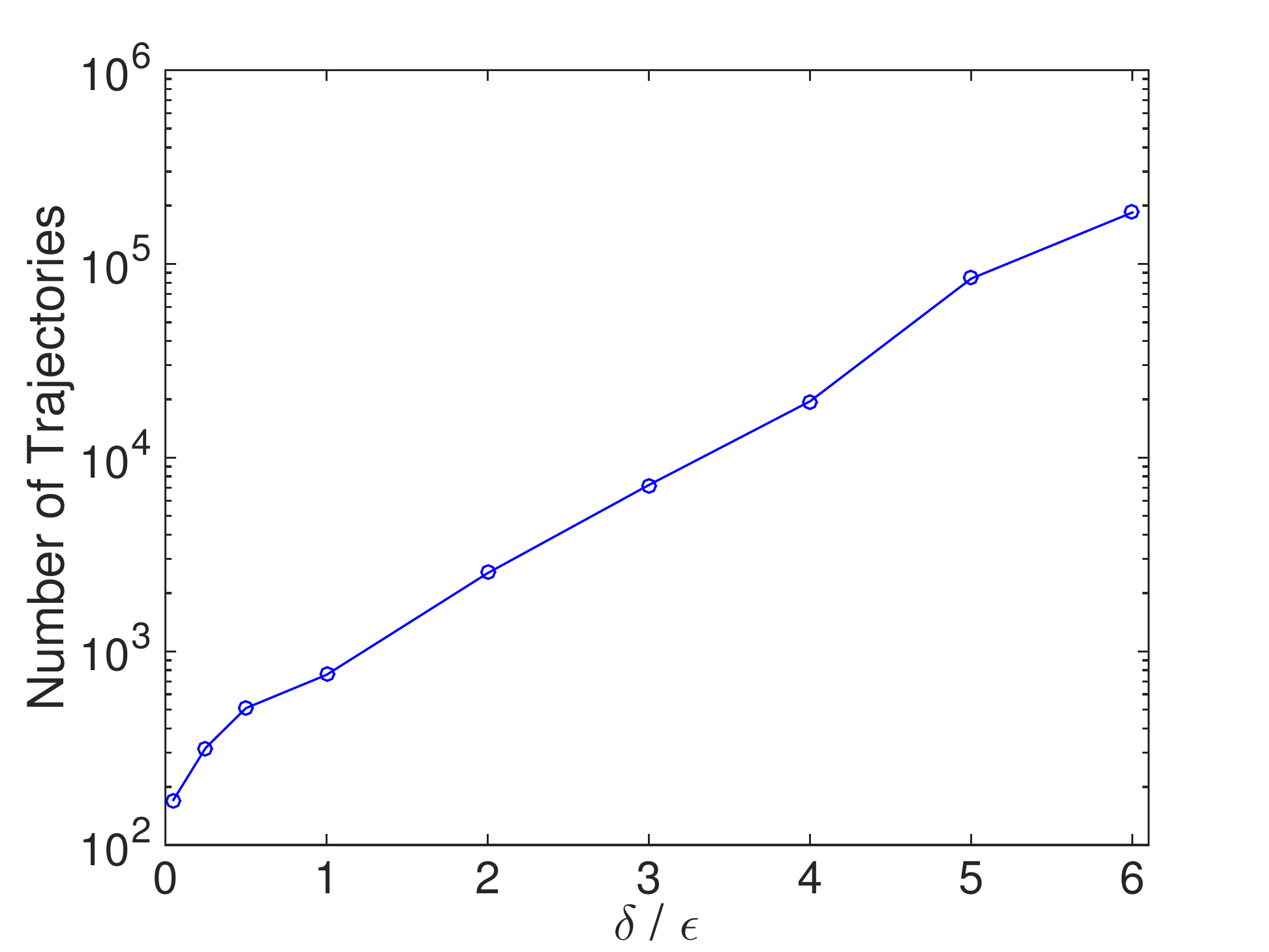}
\end{center}
\caption{Number of trajectories to achieve certain threshold with
  various $\frac{\delta}{\eps}$ in semi-log scale. Here we fix
  $\eps=0.02$, and choose $\delta= 0.001$, $0.005$, $0.01$, $0.02$,
  $0.04$, $0.06$, $0.08$, $0.1$, $0.12$. The number of trajectories
  grows exponentially with respect to $\delta$ to achieve the same
  $L^2$ relative error threshold $0.08$.}
 \label{n_traj_delta}
\end{figure}

\end{example}

\begin{example}[Avoided Crossing Regime with varying $\eps$]
  To continue, we test our algorithm in the setting of avoided
  crossing as studied in \cite{hagedorn1998landau}, where $\delta$ is
  chosen as the order of $\sqrt{\eps}$. This particular choice of
  $\delta$ generates a family of problems, where the transition rate
  gives an order one contribution due to the Landau-Zener formula.  As
  mentioned before, the diabatic representation is most useful inside
  the transition zone where the off-diagonal terms are relatively
  small, while adiabatic picture can be applied outside the transition
  zone. In this example, the propagation of the solution near the
  transition zone is considered. Since the gap is of order
  $\sqrt{\eps}$, the time for the wave packets to pass the transition zone
  is of the same order. In our numerical tests, for
  $\eps = 0.04, 0.01, 0.0025, 0.00125$, one considers
  $\delta = \sqrt{\eps}$, initial datum
  \begin{equation}
    u_{0}\left(0, x\right)=e^{-\frac{1}{2\eps}\left( x-\mu_q\right) ^{2}}e^{\frac{2i}{\eps }\left( x-\mu_q\right) },
  \end{equation}
  where $\mu_q = -2\sqrt{\eps}$, and run till the final time
  $T = 3\sqrt{\eps}$, where the wave packet passes the transition
  zone. As can be observed from Fig. \ref{avoided_stress}, the number
  of trajectories needed does not increase noticeably as $\eps$
  decreases, which shows the validity of the algorithm for small
  $\eps$.
 
  \begin{figure}[htbp] {}
    \begin{center}
      \includegraphics[scale=0.5]
      {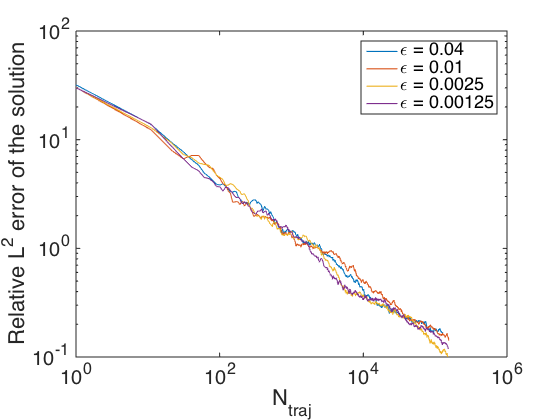}
    \end{center}
    \caption{Relative $L^2$ error of the wave function versus the
      number of trajectories for various $\eps$ in the log-log
      scale. It can be seen that as $\eps$ decreases, the number of
      trajectories does not increase noticeably for the wave packets to
      pass the transition zone.  }
 \label{avoided_stress}
\end{figure}

\end{example}

\section{Derivation of the algorithm}

\label{asym}

\subsection{A brief review of the frozen Gaussian approximation}

Recall that our algorithm is based on semiclassical approximation of
the propagation on a single diabatic surface, together with the
surface hopping to stochastically evaluate the series expansion. In
the FGA-SH method, the semiclassical approximation is based on the
frozen Gaussian approximation (\textit{aka} Herman-Kluk propagator
\cite{hermankluk}), which is a convergent semiclassical propagator of
scalar Schr\"{o}dinger equation. The Herman-Kluk propagator is
justified in \cite{kay1994integral, kay2006herman,
  swart2009mathematical} and further applied to high frequency wave
propagations \cite{fga1, fga2}.  The idea of FGA is to decompose the
initial data into a group of localized Gaussian wave packets with
fixed width, and then to propagate those Gaussians along classical
trajectories with evolving weights. To be more specific, let us recall
here the FGA integral representation for the scalar Schr\"{o}dinger
equation
\begin{equation*}
i\varepsilon \partial _{t}u\left( t,x\right) =-\frac{\varepsilon ^{2}}{2}%
\Delta _{x}u\left( t,x\right) +W\left( x\right) u\left( t,x\right) ,\
u\left( 0,x\right) =u_{\text{in}}\left( x\right), \qquad \left( t,x\right) \in 
\mathbb{R}^{+}\times \mathbb{R} ^{m}
\end{equation*}%
where $W$ is the potential (to distinguish notations from $V$ reserved for
the matrix potential). The FGA ansatz is given by%
\begin{equation*}
u_{\text{F}}\left( t,x\right) =\frac{1}{\left( 2\pi \varepsilon \right)
^{3m/2}}\int_{\mathbb{R}^{2m}}A_{t}\left(q_{0},p_{0}\right) \exp \left( \frac{i}{%
\varepsilon }\Theta_{t} \left( q_{0},p_{0},x\right) \right) \, dq_{0}dp_{0} ,
\end{equation*}%
where
\begin{equation} \label{fga_theta}
\Theta_{t} \left(q_{0},p_{0},x\right) =S_{t}\left(q_{0},p_{0}\right) +P_{t}\left(q_{0},p_{0}\right)
\left( x-Q_{t}\left(q_{0},p_{0}\right) \right) +\frac{i}{2}\lvert x-Q_{t}\left(
q_{0},p_{0}\right) \rvert ^{2}.
\end{equation}%
Notice $u_{\text{F}}$ can be viewed as a summation of localized Gaussian
wave packets. The FGA coefficients $(Q_{t},P_{t},S_{t},A_{t})$ 
satisfy the following ODEs:%
\begin{align} \label{ode_scalar}
  \frac{\rd}{\rd t}Q &=P, \\
  \frac{\rd}{\rd t}P &=-\nabla W\left( Q\right),  \\
  \frac{\rd}{\rd t}S &=\frac{1}{2}\lvert P\rvert ^{2}-W\left(
                       Q\right),  \\
  \label{ode_scalara} \frac{\rd}{\rd t}A &=\frac{1}{2}A \tr\left( Z^{-1}\bigl( \partial
                       _{z}P-i\partial _{z}Q\nabla _{Q}^{2}W\left(
                       Q\right) \bigr) \right),
\end{align}%
where we have used the shorthands $\partial _{z}=\partial _{q}-i\partial
_{p} $ and $Z=\partial _{z}\left( Q+iP\right)$. To assign the initial data of the FGA coefficients, we quote \cite[Lemma 3.1]{fga1}:
\begin{lemma} 
For $u \in L^2(\mathbb{R}^m)$, it holds
\begin{equation}
  u(x)=\frac{1}{(2 \pi \eps )^{3m/2}} \int_{\mathbb{R}^{3m}} 2^{m/2}
  u\left(y\right) \exp \left( \frac{i}{\varepsilon }\bigl( -p\left( y-q\right) +\frac{i}{2}%
    \lvert y-q\rvert ^{2}+p(x-q)+\frac{i}{2}\lvert x-q \rvert^{2} \bigr) \right)dy dq dp.
\end{equation}
\end{lemma}
As a direct consequence, the initial conditions of the action and
weight of each trajectory can be given as \eqref{initial_A}. The
$L^{2}$ error of this FGA ansatz compared to the exact solution of the
scalar Schr\"{o}dinger equation is proved to be of first order in
$\eps$ (see, e.g., \cite{swart2009mathematical}).

\subsection{A generalization of FGA to diabatic surface
  hopping} \label{asymptotic derivation} 

In this section, we ``translate'' the terms of the series
\eqref{series0}--\eqref{u_series} into a deterministic asymptotic
representation of the wave function of \eqref{schd} by using the
frozen Gaussian approximation as the semiclassical approximation.

It is clear that the first term of \eqref{u_series} needs no
``translation''. In other words,
$u^{(0)}_{0}(t,x)=\mathcal{U}_{0}(t,0)u_{0}(0,x)$ is exactly the same
as the scalar case with the potential $V_{00}(x)$, so we can directly
apply the frozen Gaussian approximation as in the previous subsection.

Now consider  
\begin{equation}
u_{0}^{(1)}(t,x)=\frac{-i \delta}{\eps}\int_{0}^{t}dt_{1}\mathcal{U}_{1}\left(
t,t_{1}\right) V_{10}(x) \mathcal{U}_{0}\left( t_{1},0\right)
u_{0}\left( 0,x\right)
\end{equation} 
and focus on the integrand at a specific $t_1$. If we view
$\mathcal{U}_{1}\left(t,t_{1}\right)$ as the propagator on the second
diabatic surface and
\begin{equation}\label{u1_initial}
V_{10}(x) \mathcal{U}_{0}\left( t_{1},0\right)u_{0}\left( 0,x\right)
\end{equation}
as the initial data at time $t_1$, the propagation from time $t_1$ to
$t$ can then be written in a FGA form
\begin{equation} \label{eq:fgaformtemp}
  \frac{1}{\left( 2\pi \varepsilon \right)
    ^{3m/2}}\int_{\mathbb{R}^{2m}}\hat A_{t} \left(q_{0},p_{0}\right) \exp \left( \frac{i}{\varepsilon } \tilde \Theta_{t}\left( q_{0},p_{0},x\right) \right) \, dq_{0}dp_{0} +O(\eps),
\end{equation}
where $ \tilde \Theta_{t}$ is defined as \eqref{fga_theta}
and the trajectory coefficients $\tilde Q_t, \tilde P_t, \tilde S_t$ and
$\hat A_t$ evlove according to the ODE system \eqref{ode_scalar}
with potential energy $W=V_{11}$ from time $t_{1}$ to $t$. To
determine the value of those at time $t_1$ (the initial time of the
propagator $\mathcal{U}(t, t_1)$), let us represent \eqref{u1_initial} in a
FGA form by evolving the data at time $t = 0$ to time $t_1$. In
particular, the last two terms $\mathcal{U}_0(t_1, 0) u_0(0, x)$ is nothing
but the scalar case with potential $W=V_{00}$, and hence can be
approximated using the scalar FGA ansatz:
\begin{equation*}
\mathcal{U}_{0}\left( t_{1},0\right)u_{0}\left( 0,x\right)= \frac{1}{\left( 2\pi \varepsilon \right)
^{3m/2}}\int_{\mathbb{R}^{2m}}A_{t_1}\left( q_{0},p_{0}\right) \exp \left( \frac{i}{%
\varepsilon }\Theta_{t_1} \left( q_{0},p_{0},x\right) \right) \, dq_{0}dp_{0} +O(\eps).
\end{equation*}
For $V_{10}(x)$, expanding around $x=Q_{t_{1}}(q_{0},p_{0})$ 
, we get
\begin{align*}
  V_{10}(x) & \mathcal{U}_{0}\left( t_{1},0\right)  u_{0}\left( 0,x\right) = \\
  = & \frac{1}{\left( 2\pi \varepsilon \right)
      ^{3m/2}}\int_{\mathbb{R}^{2m}} V_{10}\left(Q_{t_1}\left(q_0,p_0\right)\right)A_{t_1}\left( q_{0},p_0\right) \exp \left( \frac{i}{%
      \varepsilon }\Theta_{t_1} \left( q_{0},p_0,x\right) \right) , dq_{0}dp_{0}  \\
  & +\frac{1}{\left( 2\pi \varepsilon \right)
     ^{3m/2}}\int_{\mathbb{R}^{2m}}
          \left(x-Q_{t_1} \left(q_0,p_0 \right)\right) \cdot
     \nabla V_{10}\left(Q_{t_1}\left(q_0,p_0\right)\right)
     A_{t_1}\left(z_{0}\right) \exp \left( \frac{i}{%
     \varepsilon }\Theta_{t_1} \left( q_{0},p_0,x\right) \right) \, dq_{0}dp_{0}  \\
  & +\frac{1}{2\left( 2\pi \varepsilon \right)
     ^{3m/2}}\int_{\mathbb{R}^{2m}}
     \left(x-Q_{t_1} \left(q_0,p_0\right)\right) \cdot
     \nabla^2 V_{10}\left(Q_{t_1}\left(q_0,p_0\right)\right)
     \left(x-Q_{t_1}\left(q_0,p_0\right)\right) \times \\
            & \qquad \times A_{t_1}\left(q_{0},p_0\right) \exp \left( \frac{i}{%
              \varepsilon }\Theta_{t_1}\left( q_{0},p_0,x\right) \right) dq_{0}dp_{0}  +O(\eps^{2}).
\end{align*}
According to \cite[lemma 3.2]{fga1}, the second and third terms are both $O(\eps)$. Hence,
\begin{multline*}
   V_{10}(x) \mathcal{U}_{0}\left( t_{1},0\right)u_{0}\left( 0,x\right) = \\
=  \frac{1}{\left( 2\pi \varepsilon \right)
      ^{3m/2}}\int_{\mathbb{R}^{2m}} V_{10}\left(Q_{t_1}\left(q_0,p_0\right)\right)A_{t_1}\left( q_{0},p_0\right) \exp \left( \frac{i}{%
      \varepsilon }\Theta_{t_1} \left( q_{0},p_0,x\right) \right) , dq_{0}dp_{0}  +O(\eps). 
\end{multline*}
We now match the above formula with \eqref{eq:fgaformtemp} at time
$t_1$, which leads to the transition conditions of the trajectory
coefficients as
\begin{align}
  \tilde Q_{t_1} (q_0,p_0)=&Q_{t_1}(q_0,p_0), \label{q1_0}\\
  \tilde P_{t_1} (q_0,p_0)=&P_{t_1}(q_0,p_0),\label{p1_0}\\
    \tilde S_{t_1} (q_0,p_0)=&S_{t_1}(q_0,p_0),\label{s1_0}\\
  \hat A_{t_1}(q_0,p_0)=& V_{10}\left(Q_{t_1}\left(q_0,p_0\right)\right)A_{t_1}\left(q_{0},p_{0}\right) .
\end{align}
To simplify the final expression, let us define $\tilde A_{t}$ such that
\begin{equation} 
 \hat A_{t}(q_0,p_0)= V_{10}\left(Q_{t_1}\left(q_0,p_0\right)\right) \tilde A_{t}(q_0,p_0).
 \end{equation}
 Since given $q_0, p_0$ and $t_1$, $\tilde A_{t}$ only differs with
 $\hat A_{t}$ by a constant factor (independent of $t$), $\tilde A_{t}$
 follows the same ODE as $\hat A_{t}$ with initial condition
 \begin{equation}\label{a1_0}
\tilde A_{t_1}(q_0,p_0)= A_{t_1}(q_0,p_0).
 \end{equation}
 Hence, one could drop the ``$\sim$" and view $(Q_{t}, P_{t}, S_{t}, A_{t})$ as a continuous trajectory 
 at hopping time $t=t_1$. In sum, the FGA ansatz of $u_{0}^{(1)}$ reads
 \begin{multline}
 u_{\text{F}}^{(1)}\left(t,t_1,q_0,p_0\right)=  \frac{-i \delta}{\veps} \frac{1}{\left( 2\pi \varepsilon \right)
^{3m/2}}  \int_{0}^{t}dt_{1}\int_{\mathbb{R}^{2m}}dq_{0}dp_{0} 
   V_{10}\left(Q_{t_1}\left(q_0,p_0\right)\right)\times \\
   \times  A_{t}\left( q_{0},p_{0}\right)  \exp \left( \frac{i}{%
\varepsilon }\Theta_{t} \left(q_{0},p_{0},x\right) \right) \, ,
\end{multline}
where the quantities in the integrand are defined as above (solving
the ODE systems of \eqref{ode_scalar}--\eqref{ode_scalara} with the
potential $W=V_{00}$ from time $0$ to $t_1$ and with the potential
$W=V_{11}$ from time $t_1$ to $t$, with gluing conditions
\eqref{q1_0}--\eqref{s1_0} and \eqref{a1_0}.
 
Similarly, the other terms can be represented as FGA ansatz, by
applying consecutively the semiclassical approximation to each
propagators. In sum, the FGA surface hopping ansatz reads
 \begin{equation}  \label{uf}
u_{\text{F}}\left( t,x\right) =%
\begin{pmatrix}
u_{\text{F}}^{\left( 0\right) }\left( t,x\right) +u_{\text{F}}^{\left( 2\right) }\left(
t,x\right) +u_{\text{F}}^{\left( 4\right) }\left( t,x\right) +\cdots  \\ 
u_{\text{F}}^{\left( 1\right) }\left( t,x\right) +u_{\text{F}}^{\left( 3\right) }\left(
t,x\right) +u_{\text{F}}^{\left( 5\right) }\left( t,x\right) +\cdots 
\end{pmatrix}%
,
\end{equation}
with the FGA ansatz for $%
u_{\text{F}}^{\left( n\right) }$ generalized as%
\begin{align} \label{un_t}
u_{\text{F}}^{\left( n\right) }\left( t,x\right) =& \left(\frac{-i \delta}{\veps} \right)^{n} \frac{1}{\left( 2\pi \varepsilon
\right) ^{3m/2}}\int_{\mathbb{R}^{2m}} dq_{0}dp_{0}
 \int_{0<t_{1}<\cdots<t_{n}<t}dT_{n:1}
 \prod_{j=1}^{n} V_{l_{t_j}l_{t_{j-1}}} \left(Q_{t_j} \right) A_{t} \exp \left( \frac{i}{\varepsilon }%
\Theta _{t} \right) ,
 \end{align}
where $T_{n:1}:=\left( t_{n},\cdots ,t_{1}\right) $ and the phase function $\Theta_t$ is given by \eqref{theta}.

In this asymptotic derivation, the terms of order higher than $O(\eps)$ in the integrands are ignored, and hence the asymptotic error reads 
\begin{equation}
\sum_{n=0}^{\infty} \left( \frac{\delta t}{\veps} \right)^n \frac{1}{n!} \veps = \exp{ \left(\frac{\delta t}{\veps} \right)} \eps.
\end{equation}
This shows the challenge of the algorithm when $\delta$ is much larger
than $\veps$ or for very long time.

\begin{remark}
  A rigorous estimate on the asymptotic error of the surface hopping
  ansatz \eqref{uf}--\eqref{un_t} is possible following a similar
  analysis as in \cite{lu2016frozen}, we will omit the details here.
\end{remark}

\subsection{Stochastic Interpretation}

As is mentioned, the surface hopping ansatz involves integrals with respect to all of the hopping times, and hence motivates a stochastic representation of the integrals. In this section, we shall interpret the surface hopping ansartz (\ref{uf}) via the stochastic process proposed in Section \ref{sh}.

Given a time interval $\left[ 0,t\right] $ and $s$ is some time within this interval, the stochastic process is labelled as $\left(Q_{s},P_{s},l_{s}\right)$. $l_{s}$ is piecewise constant with almost surely finite many jumps. Consider a realization with
initial condition $\left(Q_{0},P_{0},l_{0}\right)=\left(q_0,p_0,0\right) $ and the number
of jumps $n$, and the discontinuity set of $l_{s}$ is denoted as $\left\{ t_{1},\cdots
t_{n}\right\} $. By the properties of the stochastic process, the
probability of no jump ($n=0$) reads%
\begin{equation}
\mathbb{P}\left( n=0\right) =\exp \left( \int_{0}^{t}\lambda _{00}\left(
Q_{s}\right) ds\right) =\exp \left( -\int_{0}^{t} \frac{\delta}{\eps} \lvert V_{01} \left(Q_{s} \right)\rvert ds\right) .
\end{equation}%
Similarly the probability with one jump ($n=1$) is given by 
\begin{align}
\mathbb{P}\left( n=1\right) =&
\int_{0}^{t}dt_{1}
\lambda_{01}(Q_{t_1})
 \exp \left(-\int_{0}^{t}  \lambda_{00}(Q_s) ds\right)  \\
=&\int_{0}^{t}dt_{1}
\frac{\delta}{\eps} \lvert V_{01} \left(Q_{ t_{1}} \right)\rvert
 \exp \left(-\int_{0}^{t} 
\frac{\delta}{\eps} \lvert V_{01} \left(Q_s \right)\rvert
ds\right) 
\end{align}%
Generally, 
\begin{equation}
P\left( n=k\right) = \int_{0<t_{1}<\cdots<t_{k}<t}dT_{k:1} \left( \frac{\delta}{\eps} \right)^k \prod_{j=1}^{k}
\lvert V_{01} \left(Q_{t_j} \right)\rvert \exp \left( -\int_{0}^{t}
\frac{\delta}{\eps} \lvert V_{01} \left(Q_{s} \right)\rvert
ds\right) , \notag
\end{equation}%
and for $k$ jumps in total, the probability density of $\left( t_{1},\cdots
,t_{n}\right) $ is given by%
\begin{equation}
\rho _{k}\left( t_{1},\cdots ,t_{k}\right) = \frac{1}{P\left( n=k\right) } \left( \frac{\delta}{ \eps} \right)^k%
\prod_{j=1}^{k}
\lvert V_{01} \left(Q_{t_j} \right)\rvert \exp \left( -\int_{0}^{t} \frac{\delta}{\eps}
\lvert V_{01} \left(Q_{s} \right)\rvert
ds\right)  \mathbb{I}_{t_{1}\leq
t_{2}\leq \cdots \leq t_{k}},
\end{equation}%
where $\mathbb{I}_{t_{1}\leq t_{2}\leq \cdots \leq t_{k}}$ is the indicator
function of the set $\left\{ \left( t_{1},\cdots ,t_{k}\right) :t_{1}\leq
t_{2}\leq \cdots \leq t_{k}\right\} $.

Now we consider a trajectory average \eqref{avg1} and denote everything
inside the square bracket of \eqref{avg1} as $\Upsilon $ , 
\begin{equation}
\Upsilon = \left( -i  \right)^n 
\prod_{j=1}^{n}\frac{V_{l_{t_j}l_{t_{j-1}}} \left(Q_{t_j} \right) }{%
\lvert V_{01} \left(Q_{t_j} \right) \rvert }%
\frac{A_{t}}{\lvert
A_{0} \rvert } \exp \left( \frac{i}{%
\varepsilon }\Theta _{t} \right) 
\exp \left( \int_{0}^{t}\frac{\delta}{\eps} \lvert
V_{\text{01}} \left(Q_s \right) \rvert ds\right) 
\begin{pmatrix}
\mathbb{I}_{n\ \text{is even}} \\ 
\mathbb{I}_{n\ \text{is odd}}%
\end{pmatrix}
\end{equation}
 and $(Q_t,P_t,l_t)$ as $\tilde z_t$ for short. We want
to show that \eqref{avg1} is a stochastic representation of the FGA ansatz \eqref{uf}.
The idea is quite straightforward. One first converts the trajectory
average into an integral with respect to $z_{0} :=(q_{0},p_0)$, which could be easily
achieved by the law of total expectation,%
\begin{align*}
\tilde{u}\left( t,x\right) &=C_{\mathcal{N}}\mathbb{E}_{\tilde z_0}\left( 
\mathbb{E}_{\tilde z_t}\left[ 
\Upsilon \mid \tilde{z}%
_{0}=\left( z_{0},0\right) \right] \right) \\
&=\frac{1}{\left( 2\pi \varepsilon \right) ^{3m/2}}\int_{\mathbb{R}%
^{2m}}dz_{0}\mathbb{E}_{\tilde{z}_{t}}\left[ \left. \lvert A_{0}(z_0) \rvert \Upsilon \right\vert \tilde{z}%
_{0}=\left( z_{0},0\right) \right] ,
\end{align*}%
where the normalization factor $C_{\mathcal{N}}$ is given as \eqref{cn}, and
the probability measure $\mathbb{P}_{0}$ of the single surface reads%
\begin{equation*}
\mathbb{P}_{0}\left( \Omega \right) =C_{\mathcal{N}}^{-1}\frac{1}{\left(
2\pi \varepsilon \right) ^{3m/2}}\int_{\Omega }dz_{0}\lvert A_{0}(z_0) \rvert .
\end{equation*}%
Next, one rewrites the condition expectation. Note that given the initial
condition, the randomness only comes from the hopping times.%
\begin{align*}
\tilde{u}\left( t,x\right) & =\frac{1}{\left( 2\pi \varepsilon \right)
^{3m/2}}\int_{\mathbb{R}^{2m}}dz_{0}\sum_{n=0}^{\infty }\mathbb{P}\left(
n\right) \int_{\left[ 0,t\right] ^{n}}dT_{n:1}\ \rho _{n}\left( t_{1},\cdots
,t_{n}\right) \lvert A_{0}\left(z_{0}\right)
\rvert \Upsilon \\
& =  \frac{1}{\left( 2\pi \varepsilon \right) ^{3m/2}}\int_{\mathbb{R}%
^{2m}}dz_{0}\sum_{n=0}^{\infty }\int_{0<t_{1}<\cdots
<t_{n}<t}dT_{n:1} \left(\frac{-i \delta}{\eps} \right)^n 
\prod_{j=1}^{n} V_{l_{t_j}l_{t_{j-1}}} \left(Q_{t_j} \right) A_{t} \exp \left( \frac{i}{\varepsilon }%
\Theta _{t} \right) 
\begin{pmatrix}
\mathbb{I}_{n\ \text{even}} \\ 
\mathbb{I}_{n\ \text{odd}}%
\end{pmatrix}
\\
& = 
\begin{pmatrix}
u_{\text{F}}^{\left( 0\right) }+u_{\text{F}}^{\left( 2\right) }+\cdots \\ 
u_{\text{F}}^{\left( 1\right) }+u_{\text{F}}^{\left( 3\right) }+\cdots%
\end{pmatrix}
 =u_{\text{F}}\left( t,x\right) .
\end{align*}%
This shows that the stochastic path integral representation can capture the
asymptotically correct solutions of the matrix Schr\"{o}dinger equation \eqref%
{schd}.

\begin{remark}
We remark that our previous transition rate \eqref{gen} depends only on the
off-diagnal term $V_{01}$. However, there are some cases when the
off-diagnal term $V_{01}$ is large while the gap between two diabatic
potential energy surfaces are even larger. Then according to previous
transition rate, we would expect excessive amount of hops happen, which,
though still gives the correct wave function, may lead to large variance and
hence low convergence rate. As one of the main issues in Monte Carlo
sampling, how to reduce the variance is always an interesting problem. But
the main purpose of this paper is to provide a general framework of the
diabatic surface hopping algorithm while efficient ways to control the
variance is beyond the scope of this paper, and will remain as further
study. Nevertheless, we remark that for the special case mentioned above,
one could choose the generator matrix as 
\begin{equation}
\lambda \left( q\right) = 
\begin{pmatrix}
\lambda _{00}\left( q\right) & \lambda _{01}\left( q\right) \\ 
\lambda _{10}\left( q\right) & \lambda _{11}\left( q\right)%
\end{pmatrix}
.
\end{equation}%
where%
\begin{equation*}
\lambda _{01}\left( q\right) =\lambda _{10}\left( q\right) =\left\{ 
\begin{array}{ll}
\frac{\delta}{\eps} \lvert \frac{ V_{01}\left( q\right) }{V_{00}\left( q\right) -V_{11}\left(
q\right) }\rvert , & \text{if }\lvert V_{00}\left( q\right)
-V_{11}\left( q\right) \rvert >1 \\ 
\frac{\delta}{\eps} \lvert V_{01}\left( q\right) \rvert , & \text{otherwise.}%
\end{array}%
\right. ,
\end{equation*}%
and%
\begin{equation*}
\lambda _{00}\left( q\right) =\lambda _{11}\left( q\right) =-\lambda
_{01}\left( q\right) .
\end{equation*}%
Note that when the gap $\lvert V_{00}\left( q\right) -V_{11}\left(
q\right) \rvert $ is large (bigger than 1), the transition rate is
modified so that the hops happen less frequently.
\end{remark}

\section{Conclusion}

A diabatic surface hopping algorithm is proposed based on time dependent
perturbation theory and semiclassical analysis. It captures the
correct asymptotic limit in the weak coupling regime. The validity of
the scheme is also verified in the avoided crossing regime for small
$\eps$ near the transition zone. The diabatic potential energy
surfaces remain smooth in the transition zone comparing to the
adiabatic energy surfaces.  It seems to us a good strategy for
application would be a mix of both adiabatic and diabatic
representations, i.e., to apply adiabatic representation away from the
transition zone and switch to diabatic picture near it.  For future
works, rigorous asymptotic analysis of the algorithm would be
interesting and better Monte Carlo sampling algorithm is also worth
exploring.

\bibliographystyle{abbrv}
\bibliography{fgash}

\end{document}